\documentclass[11pt]{article}
\usepackage{amsfonts,amsmath,amssymb}
\usepackage{mathrsfs}
\usepackage[latin1]{inputenc}
\newtheorem{theo}{Theorem}
\newenvironment{Theo}{\begin{theo}\slshape}{\end{theo}}
\newtheorem{Lemm}{Lemma}[section]

\newtheorem{Corol}{Corollary}
\newtheorem{rema}{Remark}

\newtheorem{defi}{Definition}

\newcommand{\p}{\partial}

\def\qed{\hfill$\square$\par \bigskip}
\newenvironment{Demo}[1]{{\bf Proof #1.~}}{\qed}
\newcommand{\R}{\mathbb{R}}

\newcommand{\M}{\mathcal{M}}

\newcommand{\g}{\mathrm{g}}
\newcommand{\cg}{c\mathrm{g}}
\newcommand{\dd}{\mathrm{d}}
\newcommand{\dv}{\, \mathrm{dv}_\mathrm{g}^{n}}
\newcommand{\dvc}{\, \mathrm{dv}_{c\mathrm{g}}^{n}}
\newcommand{\dvv}{\, \mathrm{dv}_\mathrm{g}^{2n-1}}
\newcommand{\ds}{\, \mathrm{d\sigma}_\mathrm{g}^{n-1}}
\newcommand{\dsc}{\, \mathrm{d\sigma}_{c\mathrm{g}}^{n-1}}
\newcommand{\dss}{\, \mathrm{d\sigma}_\mathrm{g}^{2n-2}}
\newcommand{\I}{\mathcal{I}}

\newcommand{\s}{\mathbb{S}}

\newcommand{\norm}[1]{\left\Vert#1\right\Vert}
\newcommand{\abs}[1]{\left\vert#1\right\vert}
\newcommand{\set}[1]{\left\{#1\right\}}
\newcommand{\para}[1]{\left(#1\right)}
\newcommand{\cro}[1]{\left[#1\right]}

\newcommand{\seq}[1]{\left<#1\right>}

\newcommand{\To}{\longrightarrow}
\usepackage{graphicx}
\newcommand{\dive}{\textrm{div}}

\usepackage{times}

\begin{document}
\title{Stable determination of coefficients in the dynamical anisotropic Schr\"odinger equation from the Dirichlet-to-Neumann map}
\author{\small \bf Mourad Bellassoued \\
\small University of Carthage,\\
\small Department of Mathematics,\\
\small Faculty of Sciences of Bizerte\\
\small 7021 Jarzouna Bizerte, Tunisia\\
\small mourad.bellassoued@fsb.rnu.tn
\and
\small \bf David Dos Santos Ferreira\thanks{David DSF was partially supported by ANR grant Equa-disp.} \\
\small Universit\'e Paris 13 \\
\small CNRS, UMR 7539 LAGA \\
\small 99, avenue Jean-Baptiste Cl\'ement \\
\small F-93 430 Villetaneuse, France \\
\small ddsf@math.univ-paris13.fr }
\date{}
\maketitle
\begin{abstract}
In this paper we are interested in establishing stability estimates in the inverse problem of determining on a compact Riemannian manifold
the electric potential or the conformal factor in a Schr\"odinger equation
 with Dirichlet data from measured Neumann boundary observations. This information is enclosed in the dynamical
Dirichlet-to-Neumann map associated to the
Schr\"odinger equation. We prove in dimension $n\geq 2$ that the
knowledge of the Dirichlet-to-Neumann map for the
Schr\"odinger equation uniquely determines
the electric potential and we establish H\"older-type stability  estimates in
determining the potential. We prove similar results for the determination of a conformal factor close to $1$.
\end{abstract}

\tableofcontents

\section{Introduction and main results}

This paper is devoted to the  study of the following inverse boundary value problem:  given a Riemannian manifold with boundary
determine the potential or the conformal factor of the metric in a dynamical Schr\"odinger equation  from the
observations made at the boundary. Let $(\M,\,\g)$ be a compact
Riemannian manifold with boundary $\partial \M$. All manifolds will be assumed smooth (which means $\mathcal{C}^\infty$) and oriented.
We denote by $\Delta_\g$ the Laplace-Beltrami operator associated to the metric~$\g$. In local coordinates,
$\g(x)=(\g_{jk})$, $\Delta_\g$ is given by
\begin{equation}\label{1.1}
\Delta_\g=\frac{1}{\sqrt{\det \g}}\sum_{j,k=1}^n\frac{\p}{\p
x_j}\para{\sqrt{\det \g}\,\g^{jk}\frac{\p}{\p x_k}}.
\end{equation}
Here $(\g^{jk})$ is the inverse of the metric $\g$ and $\det
\g=\det(\g_{jk})$. Let us consider the following initial boundary
value problem for the Schr\"odinger equation with bounded electric
potential $q\in L^\infty(\M)$
\begin{equation}\label{1.2}
\left\{
\begin{array}{llll}
\para{i\partial_t+\Delta_\g +q(x)}u=0,  & \textrm{in }\; (0,T)\times \M\cr
\\
u(0,\cdot )=0, & \textrm{in }\; \M\cr
\\
u=f, & \textrm{on} \,\,(0,T)\times\p \M
\end{array}
\right.
\end{equation}
where $f\in H^1((0,T)\times\p\M)$. Denote by $\nu$ the outward normal vector field along the boundary $\p \M$, so that
$\displaystyle\sum_{j,k=1}^n \g^{jk}\nu_j\nu_k=1$. Further, we may define the
dynamical Dirichlet-to-Neumann map $\Lambda_{\g,\,q}$ associated to the Schr\"odinger equation by
\begin{equation}\label{1.3}
\Lambda_{\g,\,q}f=\sum_{j,k=1}^n\nu_j\g^{jk}\frac{\p u}{\p x_k}\Big|_{(0,T)\times\p\M}.
\end{equation}
Unique determination of the metric $\g=(\g_{jk})$ from the knowledge of the Dirichlet-to-Neumann map
$\Lambda_{\g,\,q}$ is hopeless: as was noted in \cite{[SU]} in the case of the wave equation,
the Dirichlet-to-Neumann map is invariant under a gauge transformation of the metric $\g$. Namely, if one pulls back the metric $\g$ by a
diffeomorphism $\Psi:\M\to \M$ which is the identity on the boundary $\Psi|_{\p\M}={\rm Id}$ into a new metric $\Psi^* \g$, one has
$\Lambda_{\Psi^*\g,\,q}=\Lambda_{\g,\,q}$.  The inverse problem has therefore to be formulated modulo the gauge invariance. However, we will
restrict our inverse problem to a conformal class of metrics (for which there is no gauge invariance): knowing $\Lambda_{\cg,q}$, can one determine the
conformal factor $c$ and the potential $q$?
\medskip

In the case of the Schr\"odinger equation,  Avdonin and Belishev  gave an affirmative answer to this question for smooth metrics conformal to the
Euclidean metric in \cite{[AB]}. Their approach is based on the boundary control method introduced by Belishev \cite{[1]} and uses in an essential way a unique continuation property. Because of the use of this qualitative property, it seems unlikely that the boundary control method would provide
accurate stability estimates. More precisely, when $\M$ is a bounded domain of $\R^n$, and $\varrho,\,q\in\mathcal{C}^2(\overline{\M})$ are real functions, Avdonin and Belishev \cite{[AB]} show that for any fixed $T>0$ the response operator (or the Neumann-to-Dirichlet map) of the Schr\"odinger equation $(i\varrho\p_tu+\Delta u-qu)=0$ uniquely determines the coefficients $\varrho$ and $q$. The problem is reduced to recovering $\varrho,\,q$ from the boundary
spectral data. The spectral data are extracted from the response operator by the use of a variational principle.
\medskip

The uniqueness in the determination of a time-dependent electromagnetic potential, appearing in a Schr\"odinger equation in a domain with obstacles,
from the Dirichlet-to-Neumann map was proved by Eskin \cite{[Eskin]}. The main ingredient in his proof is the construction of geometrical optics solutions. In \cite{[Avdonin-al]}, Avdonin, Lenhart and Protopopescu use the so-called BC (boundary control) method to  prove that the Dirichlet-to-Neumann map determines the time-independent electrical potential in a one dimensional Schr\"odinger equation.
\medskip

The analogue problem for the wave equation has a long history. Unique determination of the metric goes back to Belishev and Kurylev \cite{[BK]} using
the boundary control method and involves works of Katchlov, Kurylev and Lassas \cite{[KKL]}, Kurylev and Lassas \cite{[KL]}) and Anderson, Katchalov,
Kurylev, Lassas and Taylor \cite{[AKKLT]}. In fact, Katchalov, Kurylev, Lassas and Mandache proved that the determination of the metric from the
Dirichlet-to-Neumann map was equivalent for the wave and Schr\"odinger equations (as well as other related inverse problems) in \cite{[KKLM]}.
Identiafiability of the potential was proved by Rakesh and Symes \cite{[Rakesh-Symes]} in the Euclidian case ($\g=e$) using complex geometrical optics solutions concentrating near lines with any direction $\omega\in \s^{n-1}$ to prove that $\Lambda_{e,q}$ determines $q(x)$ uniquely in the wave equation. This result was generalized by Ramm and Sj\"ostrand \cite{[Ramm-Sjostrand]} and Eskin \cite{[Eskin1],[Eskin2]} to the case of $q$ depending on space and time. Isakov \cite{[Isakov1]} also considered the simultaneous determination of a potential and a damping coefficient.
\medskip

As for the stability of the wave equation in the Euclidian case, we also refer to \cite{[Sun]} and \cite{[IS]}; in those papers, the Dirichlet-to-Neumann map was considered on the  whole boundary. Isakov and Sun \cite{[IS]} proved that the difference in some subdomain of two coefficients is estimated by an operator norm of the difference of the corresponding local Dirichlet-to-Neumann maps, and that
the estimate is of H\"older type. Bellassoued, Jellali and Yamamoto \cite{[Bell-Jel-Yama2]} considered the inverse problem of recovering a
time independent potential in the hyperbolic equation from the partial Dirichlet-to-Neumann map. They proved a logarithm stability estimate. Moreover in \cite{[Rakesh1]}  it is proved that if an unknown coefficient belongs to a given finite dimensional vector space, then the uniqueness
follows by a finite number of measurements on the whole boundary. In \cite{[Bellassoued-Benjoud]}, Bellassoued and Benjoud used complex geometrical
optics solutions concentring near lines in any direction to prove that the Dirichlet-to-Neumann map determines uniquely the magnetic field induced by
a magnetic potential in a magnetic wave equation.
\medskip

In the case of the anisotropic wave equation, the problem of establishing stability estimates in determining the metric was studied by Stefanov and
Uhlmann in \cite{[SU], [SU2]} for metrics close to Euclidean and generic simple metrics. In a previous paper \cite{[BellDSF]}, the authors also proved stability estimates for the wave equation in determining a conformal factor close to 1 and time independent potentials in simple geometries. We refer to this paper for a longer bibliography in the case of the wave equation.
\medskip

The inverse problem for the (dynamical) Schr\"odinger equation seems to have been a little bit less studied.
In the Euclidean case, there are extensive results by Bellassoued and Choulli \cite{[Bel-Choul]} where a Lipschitz stability estimate was proven
for time independent magnetic potentials. The stability problem in determining a time independent potential in a Schr\"odinger equation from a
single boundary measurement was studied by Baudouin and Puel~\cite{[Baudouin-Puel]}. They established Lipschitz stability estimates by a method based essentially on an appropriate Carleman inequality. In the above mentionned papers, the main assumption is that the part of the boundary where the measurement is made must satisfy a geometrical condition (related to geometric optics condition insuring observality). Recently, Bellassoued and Choulli showed in \cite{[BC]} that this geometric condition can be relaxed provided that the potential is known near the boundary. The key idea was the following : the authors used an FBI transform to change the Schr\"odinger equation near the boundary into a heat equation for which one can use a useful Carleman inequality involving a boundary term and  without any geometric condition.
\medskip

The main goal of this paper is to study the stability of the inverse problem for the dynamical anisotropic Schr\"odinger equation.
We follow the same strategy as in \cite{[BellDSF]} inspired by the works of Dos Santos Ferreira, Kenig, Salo and Uhlmann \cite{[DKSU]},
Stefanov and Uhlmann \cite{[SU],[SU2]} and Bellassoued and Choulli \cite{[Bel-Choul]}.

\subsection{Weak solutions of the Schr\"odinger equation}
First, we will consider the initial-boundary value problem for the Schr\"odinger equation on a manifold with boundary (\ref{1.2}). This initial
boundary value problem corresponds to an elliptic operator $-\Delta_\g$ given by (\ref{1.1}). We will develop an invariant approach to prove existence and uniqueness of solutions and to study their regularity proprieties.
\medskip

Before stating our first main result, we recall the following preliminaries. We refer to \cite{[Jost]} for the differential calculus of tensor fields on a Riemannian manifold. Let $(\M,\g)$ be an $n$-dimensional, $n\geq 2$, compact Riemannian manifold, with smooth boundary and smooth metric $\g$. Fix a coordinate system $x=\cro{x_1,\ldots,x_n}$ and let $\cro{\frac{\p}{\p x_1},\dots,\frac{\p}{\p x_n}}$ be the corresponding tangent vector fields. For $x\in \M$, the inner product and the norm on the tangent space $T_x\M$ are given by
\begin{gather*}
\g(X,Y)=\seq{X,Y}_\g=\sum_{j,k=1}^n\g_{jk}\alpha_j\beta_k, \\
\abs{X}_\g=\seq{X,X}_\g^{1/2},\qquad  X=\sum_{i=1}^n\alpha_i\frac{\p}{\p x_i},\quad Y=\sum_{i=1}^n
\beta_i\frac{\p}{\p x_i}.
\end{gather*}
If $f$ is a $\mathcal{C}^1$ function on $\M$, the gradient of $f$ is the vector field $\nabla_\g f$ such that
$$
X(f)=\seq{\nabla_\g f,X}_\g
$$
for all vector fields $X$ on $\M$. This reads in coordinates
\begin{equation}\label{1.4}
\nabla_\g f=\sum_{i,j=1}^n\g^{ij}\frac{\p f}{\p x_i}\frac{\p}{\p x_j}.
\end{equation}
The metric tensor $\g$ induces the Riemannian volume $\dv=\para{\det \g}^{1/2}dx_1\wedge\cdots \wedge dx_n$. We denote by $L^2(\M)$ the completion
of $\mathcal{C}^\infty(\M)$ endowed with the usual inner product
$$
\seq{f_1,f_2}=\int_\M f_1(x) \overline{f_2(x)} \dv,\qquad  f_1,f_2\in\mathcal{C}^\infty(\M).
$$
The Sobolev space $H^1(\M)$ is the completion of $\mathcal{C}^\infty(\M)$ with respect to the norm $\norm{\,\cdot\,}_{H^1(\M)}$,
$$
\norm{f}^2_{H^1(\M)}=\norm{f}^2_{L^2(\M)}+\norm{\nabla_{\g} f}^2_{L^2(\M)}.
$$
The normal derivative is
\begin{equation}\label{1.5}
\p_\nu u:=\nabla_\g u\cdot\nu=\sum_{j,k=1}^n\g^{jk}\nu_j\frac{\p u}{\p x_k}
\end{equation}
where $\nu$ is the unit outward vector field to $\p \M$.
Moreover, using covariant derivatives (see \cite{[Hebey]}), it is possible to define coordinate invariant norms in $H^k(\M)$, $k\geq 0$.
\medskip

Before stating our main results on the inverse problem, our interest will focus on the study of the initial boundary problem (\ref{1.2}), when $u$ is
a weak solution in the class $\mathcal{C}(0,T;H^1(\M))\cap\mathcal{C}^1(0,T;H^{-1}(\M))$. The following theorem gives conditions on $f$ and $q$, which
guarantee uniqueness and continuous dependence on the data of the solutions of the Schr\"odinger equation (\ref{1.2}) with non-homogenous Dirichlet
boundary condition.
\begin{Theo}\label{Th0}
Let $T>0$ be given. Suppose that $f\in H^1((0,T)\times\p\M)$ and $q\in W^{1,\infty}(\M)$. Then the unique solution $u$ of (\ref{1.2}) satisfies
\begin{align}\label{1.6}
u&\in\mathcal{C}(0,T;H^1(\M))\cap\mathcal{C}^1(0,T;H^{-1}(\M)), \\
\label{1.7}
\p_\nu u&\in L^2((0,T)\times\p\M).
\end{align}
Furthermore, there is a constant $C=C(T,\M)>0$ such that
\begin{equation}\label{1.8}
\norm{\p_\nu u}_{L^2((0,T)\times\p\M)}\leq C\norm{f}_{H^1((0,T)\times\p\M)}.
\end{equation}
The Dirichlet-to-Neumann map  $\Lambda_{\g,q}$ defined by (\ref{1.3}) is therefore continuous and we denote by
$\norm{\Lambda_{\g,q}}$ its norm in ${\cal L}(H^1((0,T)\times\p\M),L^2((0,T)\times\p\M))$.
\end{Theo}
Theorem \ref{Th0} gives a rather comprehensive treatment of the regularity problem for (\ref{1.2}) with stronger boundary condition $f$. Moreover,
our treatment clearly shows that a regularity for $f\in H^1((0,T)\times\p\M)$ is sufficient to obtain the desired interior regularity of $u$ on
$(0,T)\times\M$ while the full strength of the assumption $f\in H^1((0,T)\times\p\M)$ is used to obtain the desired boundary regularity for
$\p_\nu u$ and then the continuity of the Dirichlet-to-Neumann map $\Lambda_{\g,q}$.
\subsection{Stable determination}
In this section we state the main stability results. Let us first
introduce the admissible class of manifolds for which we can prove
uniqueness and stability results in our inverse problem. For this we
need the notion of simple manifolds \cite{[SU2]}.
\medskip

Let $(\M,\g)$ be a Riemannian manifold with boundary $\p\M$, we denote by $D$ the Levi-Civita connection on $(\M,\g)$.
For a point $x \in \p\M$, the second quadratic
form of the boundary
$$
\Pi(\theta,\theta)=\seq{D_\theta\nu,\theta}_\g,\quad \theta\in T_x(\p\M)
$$
is defined on the space $T_x(\p\M)$. We say that the boundary is strictly convex if the form is positive-definite for all $x \in \p\M$.
\begin{defi}
We say that the Riemannian manifold $(\M,\g)$ (or that the metric $\g$) is simple in $\M$, if $\p \M$ is
strictly convex with respect to $\g$, and for any $x\in \M$, the exponential map
$\exp_x:\exp_x^{-1}(\M)\To \M$ is a diffeomorphism. The latter means that every two points $x; y \in \M$ are joined by a unique geodesic smoothly
depending on $x$ and $y$.
\end{defi}
Note that if $(\M,\g)$ is simple, one can extend it to a simple manifold $\M_{1}$ such that $\M_1\supset\overline{\M}$.
\medskip

Let us now introduce the admissible set of potentials $q$ and the admissible set of conformal factors $c$. Let $M_0>0$, $k\geq 1$ and
$\varepsilon>0$ be given, set
\begin{equation}\label{1.11}
\mathscr{Q}(M_0)=\set{q\in W^{1,\infty}(\M),\,\norm{q}_{W^{1,\infty}(\M)}\leq M_0},
\end{equation}
and
\begin{multline}\label{1.12}
\mathscr{C}(M_0,k,\varepsilon)=\\
 \set{c\in\mathcal{C}^\infty(\M),\,\,c>0\,\,\textrm{in}\,\overline{\M},\,\,
\norm{1-c}_{\mathcal{C}^1(\M)}\leq\varepsilon,\,\,\norm{c}_{\mathcal{C}^k(\M)}\leq M_0}.
\end{multline}
The main results of this paper are as follows.
\begin{Theo}\label{Th1}
Let $(\M,\g)$ be a simple compact Riemannian manifold with boundary of dimension $n \geq 2$ and let $T>0$. There exist constants
$C > 0$ and $s\in (0,1)$ such that for any $q_1,\,q_2\in\mathscr{Q}(M_0)$, $q_1=q_2$ on $\p\M$, we have
\begin{equation}\label{1.13}
\norm{q_1-q_2}_{L^2(\M)}\leq C
\norm{\Lambda_{\g,q_1}-\Lambda_{\g,q_2}}^{s}
\end{equation}
where $C$ depends on $\M$, $T$, $M_0$, $n$, $\alpha$ and $s$.
\end{Theo}
By Theorem \ref{Th1}, we can readily derive the following uniqueness result
\begin{Corol}
Assume that $T>0$. Let $q_1,\,q_2\in\mathscr{Q}(M_0)$, $q_1=q_2$ on $\p\M$. Then $\Lambda_{\g,q_1}=\Lambda_{\g,q_2}$ implies
$q_1 = q_2$ everywhere in $\M$.
\end{Corol}
\begin{Theo}\label{Th2}
Let $(\M,\g)$ be a simple compact Riemannian manifold with boundary of dimension $n \geq 2$ and let $T>0$.
There exist $k\geq 1$, $\varepsilon>0$, $0<s<1$ and $C>0$ such that for any $c\in\mathscr{C}(M_0,k,\varepsilon)$ with $c=1$ near
the boundary $\p\M$, the following estimate holds true
\begin{equation}\label{1.14}
\norm{1-c}_{L^2(\M)}\leq C\norm{\Lambda_{\g}-\Lambda_{\cg}}^{s}
\end{equation}
where $C$ depends on $\M$, $M_0$, $n$,
$\varepsilon$, $k$ and $s$.
\end{Theo}
By Theorem \ref{Th2}, we can readily derive the following uniqueness result
\begin{Corol}
Let $(\M,\g)$ be a simple compact Riemannian manifold with boundary of dimension $n \geq 2$ and let $T>0$.
There exist $k\geq 1$, $\varepsilon>0$, such that for any $c\in\mathscr{C}(M_0,k,\varepsilon)$ with $c=1$ near
the boundary $\p\M$, we have that $\Lambda_{\cg}=\Lambda_{\g}$ implies
$c =1$ everywhere in $\M$.
\end{Corol}
\medskip

Our proof is inspired by techniques used by Stefanov and Uhlmann \cite{[SU2]},  and Dos Santos Ferreira-Kenig-Salo-Uhlmann \cite{[DKSU]}
which prove uniqueness theorems for an inverse problem related to an elliptic equation. Their idea in turn goes back to the pioneering work of Calder\'on \cite{[Calderon]}. We also refer to Bukhgeim and Uhlamnn \cite{[Bukhgeim-Uhlmann]}, Cheng and Yamamoto \cite{[CY]}, Hech-Wang
\cite{[Hech-Wang]} and Uhlmann \cite{[Uhlmann]} as a survey.
\medskip

The outline of the paper is as follows. In section 2  we collect some of the formulas needed in the paper. In section 3 we study the Cauchy problem
for the Schr\"odinger equation and we prove Theorem \ref{Th0}. In section 4 we construct special geometrical optic solutions to Schr\"odinger equations.
In section 5 and 6, we establish stability estimates for related integrals over geodesics crossing $\M$ and prove our main results.
\section{Geodesical ray transform on a simple manifold}
\setcounter{equation}{0}
In this section we first collect some  formulas needed in the rest of this paper and introduce the geodesical ray transform. Denote by $\dive X$ the divergence of a vector field $X\in H^1(T\M)$ on $\M$, i.e. in local coordinates,
\begin{equation}\label{2.1}
\dive X=\frac{1}{\sqrt{\det\g}}\sum_{i=1}^n\p_i\para{\sqrt{\det\g}\,\alpha_i},\quad X=\sum_{i=1}^n\alpha_i\frac{\p}{\p x_i}.
\end{equation}
If $X\in H^1(T\M)$ the divergence formula reads
\begin{equation}\label{2.2}
\int_\M\dive X \dv=\int_{\p \M}\seq{X,\nu}\ds
\end{equation}
and for $f\in H^1(\M)$ Green's formula reads
\begin{equation}\label{2.3}
\int_\M\dive X\,f\dv=-\int_\M\seq{X,\nabla_\g f}_\g\dv+\int_{\p \M}\seq{X,\nu} f\ds.
\end{equation}
Then if $f\in H^1(\M)$ and $w\in H^2(\M)$, the following identity holds
\begin{equation}\label{2.4}
\int_\M\Delta_\g w f\dv=-\int_\M\seq{\nabla_\g w,\nabla_\g f}_\g\dv+\int_{\p \M}\p_\nu w f \ds.
\end{equation}
Let $v\in\mathcal{C}^1(\M)$ and $N$ be a smooth real vector field. The following identity holds true (see \cite{[Yao]})
\begin{multline}\label{2.4'}
\big\langle\nabla_\g v, \nabla_\g \big(\seq{N,\nabla_\g\overline{v}}_\g\big)\big\rangle_\g  \\ 
= DN(\nabla_\g v, \nabla_\g\overline{v})+
\frac{1}{2}\dive\big(\abs{\nabla_\g v}_\g^2N\big) -\frac{1}{2}\abs{\nabla_\g v}_\g^2\dive(N)
\end{multline}
where $D$ is the Levi-Civita connection.\\
For $x\in \M$ and $\theta\in T_x\M$ we denote by $\gamma_{x,\theta}$ the unique geodesic starting at the point $x$ in the direction $\theta$.  We
consider
\begin{align*}
S\M&=\set{(x,\theta)\in T\M;\,\abs{\theta}_\g=1}, \\
S^*\M&=\set{(x,p)\in T^*\M;\,\abs{p}_\g=1}
\end{align*}
the sphere bundle and co-sphere bundle of $\M$. The exponential map $\exp_x:T_x\M\To \M$ is given by
\begin{equation}\label{2.5}
\exp_x(v)=\gamma_{x,v}(\abs{v}_\g v)=\gamma_{x,v}(rv),\quad r=\abs{v}_\g.
\end{equation}
A compact Riemannian manifold $(\M,\, \g)$ with boundary is called a convex non-trapping
manifold, if it satisfies two conditions:
\begin{enumerate}
    \item[(a)] the boundary $\p \M$ is strictly convex, i.e., the second fundamental form of the boundary is positive definite at every
    boundary point,
    \item[(b)] for every point $x \in \M$ and every vector $\theta\in T_x\M$, $\theta\neq 0$, the maximal
    geodesic $\gamma_{x,\theta}(t)$ satisfying the initial conditions $\gamma_{x,\theta}(0) = x$ and $\dot{\gamma}_{x,\theta}(0) = \theta$ is defined on
    a finite segment $[\tau_{-}(x,\theta), \tau_{+}(x,\theta)]$. We recall that a geodesic $\gamma: [a, b] \To M$ is maximal
    if it cannot be extended to a segment $[a-\varepsilon_1, b+\varepsilon_2]$, where $\varepsilon_i \geq 0$ and $\varepsilon_1 + \varepsilon_2 > 0$.
\end{enumerate}
The second condition is equivalent to all geodesics having finite length in $\M$.
\medskip

An important subclass of convex non-trapping manifolds are simple manifolds. We way that a compact Riemannian
manifold $(\M, \g)$ is simple if it satisfies the following properties
\begin{enumerate}
     \item[(a)] the boundary is strictly convex,
     \item[(b)] there are no conjugate points on any geodesic.
\end{enumerate}
A simple $n$-dimensional Riemannian manifold is diffeomorphic to a closed ball in $\R^n$, and any pair of
points in the manifold are joined by an unique geodesic.
\medskip

Now, we introduce the submanifolds of inner and outer vectors of $S\M$
\begin{equation}\label{3.1}
\p_{\pm}S\M =\set{(x,\theta)\in S\M,\, x \in \p \M,\, \pm\seq{\theta,\nu(x)}< 0}
\end{equation}
where $\nu$ is the unit outer normal to the boundary. Note that $\p_+ S\M$ and $\p_-S\M$ are
compact manifolds with the same boundary $S(\p \M)$, and $\p S\M = \p_+ S\M\cup \p_- S\M$.
For $(x,\theta)\in\p_+ S\M$,  we denote by $\gamma_{x,\theta} : [0,\tau_+(x,\theta)] \To \M$ the maximal
geodesic satisfying the initial conditions $\gamma_{x,\theta}(0) = x$ and $\dot{\gamma}_{x,\theta}(0) = \theta$.
Let $\mathcal{C}^\infty(\p_+ S\M)$ be the space of smooth functions on the manifold $\p_+S\M$. The ray
transform (also called geodesic X-ray transform) on a convex non trapping manifold $\M$ is the linear operator
\begin{equation}\label{3.2}
\I:\mathcal{C}^\infty(\M)\To \mathcal{C}^\infty(\p_+S\M)
\end{equation}
defined by the equality
\begin{equation}\label{3.3}
\I f(x,\theta)=\int_0^{\tau_+(x,\theta)}f(\gamma_{x,\theta}(t))\, \dd t.
\end{equation}
The right-hand side of (\ref{3.3}) is a
smooth function on $\p_+S\M$ because the integration limit $\tau_+(x,\theta)$ is a smooth function on $\p_+S\M$, see Lemma 4.1.1 of \cite{[Sh]}.
The ray transform on a convex non trapping manifold $\M$ can be extended as a bounded operator
\begin{equation}\label{3.4}
\I:H^k(\M)\To H^k(\p_+S\M)
\end{equation}
for every integer $k\geq 1$, see Theorem 4.2.1 of \cite{[Sh]}.\\

The Riemannian scalar product on $T_x\M$ induces the volume form on $S_x\M$,
denoted by $d\omega_x(\theta)$ and given by
$$
\dd \omega_x(\theta)=\sum_{k=1}^n(-1)^k\theta^k \dd \theta^1\wedge\cdots\wedge \widehat{\dd \theta^k}\wedge\cdots\wedge \dd \theta^n.
$$
We introduce the volume form $\dvv$ on the manifold $S\M$ by
$$
\dvv (x,\theta)=\abs{\dd\omega_x(\theta)\wedge \dv}
$$
where $\dv$ is the Riemannnian volume form on $\M$. By Liouville's theorem, the form $\dvv$ is preserved by the geodesic flow. The
corresponding volume form on the boundary $\p S\M =\set{(x,\theta)\in S\M,\, x\in\p \M}$ is given
by
$$
\dss=\abs{\dd\omega_x(\theta)\wedge \ds}
$$
where $\ds$ is the volume form of $\p \M$.
\medskip

Let $L^2_\mu(\p_+S\M)$ be the space of square integrable functions with respect to the measure $\mu(x,\theta)\dss$ with
$\mu(x,\theta)=\abs{\seq{\theta,\nu(x)}}$. This real Hilbert space is endowed with the scalar product
\begin{equation}\label{3.5}
\seq{u,v}_{L^2_\mu(\p_+S\M)}=\int_{\p_+S\M}u(x,\theta)v(x,\theta)\mu(x,\theta)\dss.
\end{equation}
The ray transform $\I$ is a bounded operator from $L^2(\M)$ into
$L^2_\mu(\p_+S\M)$. The adjoint $\I^*:L^2_\mu(\p_+S\M)\To L^2(\M)$ is given by
\begin{equation}\label{3.6}
\I^*\psi(x)=\int_{S_x\M}\psi^*(x,\theta)\, \dd\omega_x(\theta)
\end{equation}
where $\psi^*$ is the extension of the function $\psi$ from $\p_+S\M$ to $S\M$ constant on every orbit of the geodesic flow, i.e.
$$
\psi^*(x,\theta)=\psi(\gamma_{x,\theta}(\tau_+(x,\theta))).
$$
Let $(\M,\g)$ be a simple metric, we assume that $\g$ extends smoothly as a simple metric on $\M_1\Supset \M$. Then there exist
$C_1>0, C_2>0$ such that
\begin{equation}\label{3.7}
C_1\norm{f}_{L^2(\M)}\leq\norm{\I^*\I(f)}_{H^1(\M_1)}\leq C_2\norm{f}_{L^2(\M)}
\end{equation}
for any $f\in L^2(\M)$.If $V$ is an open set of the simple Riemannian manifold $(\M_{1},\g)$, the normal operator $\I^*\I$ is an
elliptic  pseudodifferential operator of order $-1$ on $V$ whose principal symbol is a multiple of $\abs{\xi}_{\g}$ (see \cite{[SU2]}).
Therefore there exists a constant $C_k>0$ such that for all $f\in H^k(V)$ compactly supported in $V$
\begin{equation}\label{3.8}
\norm{\I^*\I(f)}_{H^{k+1}(\M_{1})}\leq C_k\norm{f}_{H^k(V)}.
\end{equation}
\section{The Cauchy problem for the Schr\"odinger equation}
\setcounter{equation}{0}
In this section we will establish existence, uniqueness and continuous dependence on the data of solutions to the Schr\"odinger equation (\ref{1.2}) with non-homogenous Dirichlet boundary condition $f\in H^1((0,T)\times\p\M)$. We will use the method of transposition, or adjoint
isomorphism of equations, and we shall solve the case of non-homogenous Dirichlet boundary conditions under stronger assumptions on the data
than those in \cite{[Baudouin-Puel]}.
\medskip

Let us first review the classical well-posedness results for the Schr\"odinger equation with homogenous boundary conditions. After applying the
transposition method, we establish Theorem \ref{Th0}.
\subsection{Homogenous boundary condition}
Let us consider the following initial and homogenous boundary value problem for the Schr\"odinger equation
\begin{equation}\label{7.1}
\left\{
\begin{array}{llll}
\para{i\partial_t+\Delta_\g+q(x)}v(t,x)=F(t,x)  & \textrm{in }\,\,(0,T)\times\M,\cr
\\
v(0,x)=0 & \textrm{in }\,\,\M,\cr
\\
v(t,x)=0 & \textrm{on } \,\, (0,T)\times\p\M.
\end{array}
\right.
\end{equation}
Firstly, it is well known that if $F\in L^1(0,T;L^2(\M))$ then (\ref{7.1}) admits an unique weak solution
\begin{equation}\label{7.2}
v\in \mathcal{C}\para{0,T;L^2(\M)}.
\end{equation}
If we multiply both sides of the first equation (\ref{7.1}) by $\overline{v}$ and integrate over $\M$, we obtain
\begin{multline}\label{7.3}
\textrm{Im}\cro{\int_\M
i\p_tv(t)\overline{v}\dv-\int_\M\abs{\nabla_\g v(t,x)}_g^2+q\abs{v(t,x)}^2\dv} \\
=\frac{1}{2}\frac{\dd}{\dd t}\int_\M\abs{v(t,x)}^2\dv=\textrm{Im}\int_\M \para{F(t,x)\overline{v}(t,x)+q\abs{v(t,x)}^2}\dv.
\end{multline}
take $\alpha_0(t)=\norm{v(t)}_{L^2(\M)}$, $t\in(0,T)$, we get
\begin{equation}\label{7.4}
\frac{\dd}{\dd t}\para{\alpha_0^2(t)}\leq
C\para{\norm{F(t,\cdot)}_{L^2(\M)}\alpha_0(t)+\alpha_0^2(t)},\quad\forall t\in(0,T),
\end{equation}
which implies that $\alpha'_0(t)\leq C\para{\norm{F(t,\cdot)}_{L^2(\M)}+\alpha_0(t)}$ and
\begin{equation}\label{7.5}
\norm{v(t)}_{L^2(\M)}\leq
C_T\int_0^T\norm{F(t,\cdot)}_{L^2(\M)} \, \dd t,\quad\forall t\in(0,T).
\end{equation}

\begin{Lemm}\label{L7.1}
Let $T>0$ and $q\in W^{1,\infty}(\M)$. Suppose that $F\in L^1(0,T;H^1_0(\M))$. Then the unique solution $v$ of \eqref{7.1} satisfies
\begin{equation}\label{7.6}
v\in \mathcal{C}(0,T;H^1_0(\M)).
\end{equation}
Furthermore there is a constant $C>0$ such that for any $F\in L^1(0,T;H^1_0(\M))$, we have
\begin{equation}\label{7.7}
\norm{v(t,\cdot)}_{H^1_0(\M)}\leq C\norm{F}_{L^1(0,T;H^1_0(\M))}.
\end{equation}
\end{Lemm}
\begin{Demo}{}
Using the classical result of existence and uniqueness of weak solutions in Cazenave and Haraux \cite{[Cazenave-Haraux]} (set for abstract evolution
equations), we obtain
\begin{equation}\label{7.8}
v\in \mathcal{C}(0,T;H^1_0(\M)).
\end{equation}
Multiplying the first equation of (\ref{7.1}) by $\Delta_\g\overline{v}$ and using Green's formula, we get
\begin{align}\label{7.9}
\textrm{Im}\bigg[\int_\M i\p_tv(t)&\Delta_\g\overline{v}\dv-\int_\M\abs{\Delta_\g v(t)}^2+qv\Delta_\g\overline{v}\dv\bigg] \\ \nonumber
&=-\frac{1}{2}\frac{d}{\dd t}\int_\M\abs{\nabla_\g v(t)}^2\dv+\textrm{Im}\int_\M \seq{\nabla_\g(qv),\nabla_\g \overline{v}}_\g\dv \\ \nonumber
&=\textrm{Im}\int_\M \seq{\nabla_\g F(t,x),\nabla_\g\overline{v(t)}}\dv.
\end{align}
Let $\alpha_1(t)=\norm{\nabla_\g v(t)}_{L^2(\M)}$, $t\in(0,T)$. Then, by (\ref{7.9}), we have
\begin{equation}\label{7.10}
\frac{\dd}{\dd t}\para{\alpha_1^2(t)}\leq
C\para{\norm{F(t,\cdot)}_{H^1_0(\M)}\alpha_1(t)+\alpha_1^2(t)},\quad\forall t\in(0,T),
\end{equation}
which implies that $\alpha'_1(t)\leq C\para{\norm{F(t,\cdot)}_{H^1_0(\M)}+\alpha_1(t)}$ and by Gronwall's lemma we find
\begin{equation}\label{7.11}
\norm{v(t)}_{H^1_0(\M)}\leq
C_T\int_0^T\norm{F(t,\cdot)}_{H^1_0(\M)}\, \dd t,\quad\forall t\in(0,T).
\end{equation}
The proof of (\ref{7.7}) is complete.
\end{Demo}
\begin{Lemm}\label{L7.2}
Let $T>0$ and $q\in L^{\infty}(\M)$. Suppose that $F\in W^{1,1}(0,T;L^2(\M))$ such that $F(0,\cdot)\equiv0$. Then the unique solution $v$ of \eqref{7.1}
satisfies
\begin{equation}\label{7.12}
v\in \mathcal{C}^1(0,T;L^2(\M))\cap\mathcal{C}(0,T;H^2(\M)\cap H^1_0(\M)).
\end{equation}
Furthermore there is a constant $C>0$ such that for any $\eta>0$ small, we have
\begin{equation}\label{7.13}
\norm{v(t,\cdot)}_{H^1_0(\M)}\leq C\para{\eta\norm{\p_tF}_{L^1(0,T;L^2(\M))}+\eta^{-1}\norm{F}_{L^1(0,T;L^2(\M))}}.
\end{equation}
\end{Lemm}
\begin{Demo}{}
If we consider the equation satisfied by $\p_t v$, (\ref{7.2}) provides the following regularity
   $$v\in {\cal C}^{1}(0,T;L^2(\M)).$$
Furthermore by (\ref{7.5}), there is a constant $C>0$ such that the following estimate holds true
\begin{equation}\label{7.14}
\norm{\p_t v(t,\cdot)}_{L^2(\M)}\leq
C\int_0^T\norm{\p_tF(t,\cdot)}_{L^2(\M)} \, \dd t\qquad\forall t\in(0,T).
\end{equation}
Then, by (\ref{7.1}), we see that  $\Delta_\g v=-i\p_t v+F \in {\cal C}(0,T;L^2(\M))$ and therefore
$$v\in \mathcal{C}(0,T;H^2(\M)).$$
This complete the proof of (\ref{7.12}).\\
Next, multiplying the first equation of (\ref{7.1}) by $\overline{v}$ and integrating by parts, we obtain
\begin{align}\label{7.15}
\textrm{Re}\bigg[\int_\M i\p_tv(t,x)\overline{v}(t,x)&\dv-\int_\M\!\!\para{\abs{\nabla_\g v(t)}^2-q\abs{v}^2}\dv\bigg] \cr
&= \textrm{Re}\!\int_\M F(t,x)\overline{v}(t,x)\dv\cr
&= \textrm{Re}\!\int_\M\!\!\para{\int_0^t\!\p_tF(s,x) \, \dd s}\overline{v}(t,x)\dv.
\end{align}
Take now $\alpha_1(t)=\norm{\nabla_\g v(t)}_{L^2(\M)}$. Then there exists a constant $C>0$ such that the following estimate holds true
\begin{align}\label{7.16}
\alpha^2_1(t)&\leq
C\Bigg[\norm{\p_tv(t,\cdot)}_{L^2(\M)}\norm{v(t,\cdot)}_{L^2(\M)}+\norm{v(t)}^2_{L^2(\M)}\cr
&\quad +\int_0^T\!\!\!\int_\M\abs{v(t,x)\p_tF(s,x)}\dv \, \dd s\Bigg].
\end{align}
Using (\ref{7.14}) and (\ref{7.5}), we get
\begin{equation}\label{7.17}
\alpha^2_1(t)\leq C\cro{\norm{\p_tF}_{L^1(0,T;L^2(\M))}\norm{F}_{L^1(0,T;L^2(\M))}+\norm{F}_{L^1(0,T;L^2(\M))}^2}.
\end{equation}
Thus, we deduce (\ref{7.13}), and this concludes the proof of Lemma \ref{L7.2}.
\end{Demo}
\begin{Lemm}\label{L7.3}
Let $T>0$, $q\in W^{1,\infty}(\M)$ be given and let $\mathscr{H}=L^1(0,T;H^1_0(\M))$ or $\mathscr{H}=H^1_0(0,T;L^2(\M))$. Then the mapping
$F\mapsto\p_\nu v$ where $v$ is the unique  solution to \eqref{7.1} is linear and continuous from $\mathscr{H}$ to $L^2((0,T)\times\p\M)$. Furthermore,
there is a constant $C>0$ such that
\begin{equation}\label{7.18}
\norm{\p_\nu v}_{L^2((0,T)\times\p\M)}\leq C\norm{F}_{\mathscr{H}}.
\end{equation}
\end{Lemm}
\begin{Demo}{}
Let $N$ be a $\mathcal{C}^2$ vector field on $\overline{\M}$ such that
\begin{equation}\label{7.19}
N(x)=\nu(x),\quad x\in\p\M;\qquad \abs{N(x)}_\g\leq 1,\quad x\in\M.
\end{equation}
Multiply both sides of the first equation in (\ref{7.1}) by $\seq{N,\nabla_\g\overline{v}}_\g$ and integrate over $(0,T)\times\M$, this gives
\begin{align}\label{7.20}
\int_0^T\!\!\!\int_\M F(t,x)&\seq{N,\nabla_\g\overline{v}}_\g \dv \, \dd t \cr
&=\int_0^T\!\!\!\int_\M i\p_tv\seq{N,\nabla_\g\overline{v}}_\g\dv \, \dd t
+\int_0^T\!\!\!\int_\M\Delta_\g\seq{N,\nabla_\g\overline{v}}_\g\dv \, \dd t \cr &\quad +\int_0^T\!\!\!\int_\M q(x)v\seq{N,\nabla_\g\overline{v}}_\g\dv
\, \dd t\cr
&=I_1+I_2+I_3.
\end{align}
Consider the first term on left side of (\ref{7.20}); integrating by parts with respect $t$, we get
\begin{align}\label{7.21}
I_1&=\int_0^T\!\!\!\int_\M i\p_tv\seq{N,\nabla_\g\overline{v}}_\g\dv \, \dd t \displaybreak[1] \\ \nonumber
&=i\cro{\int_\M v\seq{N,\nabla_\g\overline{v}}_\g\dv}_0^T-i\int_0^T\!\!\!\int_\M v\seq{N,\nabla_\g\p_t\overline{v}}_\g\dv \, \dd t \displaybreak[1] \\ \nonumber
&=i\int_\M v(T,x)\seq{N,\nabla_\g\overline{v}(T,x)}_\g\dv
-i\int_0^T\!\!\!\int_\M\seq{N,\nabla_\g(v\p_t\overline{v})}_\g\dv \, \dd t \cr &\quad +i\int_0^T\!\!\!\int_\M\p_t\overline{v}\seq{N,\nabla_\g v}_\g
\dv\, \dd t.
\end{align}
Then, by (\ref{2.3}), we obtain
\begin{align*}
\textrm{Re}&\bigg[\int_0^T\!\!\!\int_\M i\p_tv\seq{N,\nabla_\g\overline{v}}_\g\dv \, \dd t\bigg]  \displaybreak[1]  \\
&=i\int_\M v(T,x)\seq{N,\nabla_\g\overline{v}(T,x)}_\g\dv +i\int_0^T\!\!\!\int_\M\dive(N) v\p_t\overline{v}\dv \, \dd t  \displaybreak[1]  \\
&\quad -i\cro{\int_0^T\!\!\!\int_{\p\M} v\p_t\overline{v}\ds \, \dd t} \displaybreak[1] \\
&=i\int_\M \!v(T,x)\seq{N,\nabla_\g\overline{v}(T,x)}_\g\dv +\int_0^T\!\!\!\int_\M\!\seq{\nabla_\g\overline{v},\nabla_\g\para{\dive(N)v}}_\g\dv \,
\dd t  \displaybreak[1] \\
&\quad+\int_0^T\!\!\!\int_\M F\dive(N)v\dv \, \dd t -\int_0^T\int_\M q\dive(N)\abs{v}^2\dv \, \dd t  \displaybreak[1] \\
&\quad-\Bigg[i\int_0^T\!\!\!\int_{\p\M} v\p_t\overline{v}\ds \, \dd t +\int_0^T\!\!\!\int_{\p\M} \p_\nu\overline{v}v\dive(N)\ds \, \dd t \Bigg].
\end{align*}
The last term vanishes, using (\ref{7.13}) or (\ref{7.7}), we conclude that
\begin{equation}\label{7.23}
\abs{\textrm{Re}\,I_1}\leq C\norm{F}^2_{\mathscr{H}}.
\end{equation}
On the other hand, by Green's theorem, we get
\begin{align*}
I_2&=\int_0^T\!\!\!\int_\M\Delta_\g v\seq{N,\nabla_\g\overline{v}}_\g\dv \, \dd t \cr
&=-\int_0^T\!\!\!\int_\M\seq{\nabla_\g v, \nabla_\g (\seq{N,\nabla_\g\overline{v}}_\g)}_\g\dv \, \dd t
+\int_0^T\!\!\!\int_{\p\M}\abs{\p_\nu v}^2\ds \, \dd t.
\end{align*}
Thus by (\ref{2.4'}), we deduce
\begin{align*}
I_2 &=\int_0^T\!\!\!\int_{\p\M}\abs{\p_\nu v}^2\ds \, \dd t-\frac{1}{2}\int_0^T\!\!\!\int_{\p\M}\abs{\nabla_g v}^2\ds \, \dd t\cr
&\quad +\int_0^T\!\!\!\int_\M D_{\g}N(\nabla_\g v, \nabla_\g\overline{v})\dv \, \dd t
-\frac{1}{2}\int_0^T\!\!\!\int_\M\abs{\nabla_\g v}_\g^2\dive(N)\dv \, \dd t.
\end{align*}
Using the fact
$$
\abs{\nabla_\g v}_\g^2=\abs{\p_\nu v}^2+\abs{\nabla_\tau v}^2_\g=\abs{\p_\nu v}^2,\quad x\in\p\M
$$
where $\nabla_\tau$ is the tangential gradient on $\p\M$, we get
\begin{align}\label{7.26}
\textrm{Re}\,I_2&=\frac{1}{2}\int_0^T\!\!\!\int_{\p\M}\abs{\p_\nu v}^2\ds \, \dd t +\int_0^T\!\!\!\int_\M D_{\g}N(\nabla_\g v, \nabla_\g\overline{v})\dv
\, \dd t\cr
&\quad -\frac{1}{2}\int_0^T\!\!\!\int_\M\abs{\nabla_\g v}_\g^2\dive(N)\dv \, \dd t.
\end{align}
Finally by Lemma \ref{L7.1} and \ref{L7.2}, we have
\begin{equation}\label{7.27}
\abs{\textrm{Re}\, I_3}\leq \norm{F}^2_{\mathscr{H}}.
\end{equation}
Collecting (\ref{7.27}), (\ref{7.26}), (\ref{7.23}) and (\ref{7.20}), we obtain
\begin{equation}\label{7.28}
\int_0^T\!\!\!\int_{\p\M}\abs{\p_\nu v}^2\ds dt\leq C\norm{F}^2_{\mathscr{H}}.
\end{equation}
This completes the proof of (\ref{7.18}).
\end{Demo}
\subsection{Non-homogenous boundary condition}
We now turn to the non-homogenous Schr\"odinger problem (\ref{1.2}). Let
$$ 
\mathscr{H}=L^1(0,T;H^1_0(\M)) \text{ or } \mathscr{H}=H^1_0(0,T;L^2(\M)). 
$$
By $\para{\cdot,\cdot}_{\mathscr{H}',\mathscr{H}}$, we denote the dual pairing between $\mathscr{H}'$ and $\mathscr{H}$.
\begin{defi}
Let $T>0$, $q\in W^{1,\infty}(\M)$ and $f\in L^2((0,T)\times\p\M)$, we say that $u\in\mathscr{H}'$
is a solution of \eqref{1.2} in the transposition sense if for any $F\in \mathscr{H}$ we have
\begin{equation}\label{7.29}
\para{u,\,F}_{\mathscr{H}',\mathscr{H}}=\int_0^T\!\!\!\int_{\p\M}
f(t,x)\p_\nu\overline{v}(t,x)\ds \, \dd t
\end{equation}
where $v=v(t,x)$ is the solution of the homogenous boundary value problem
\begin{equation}\label{7.33}
\left\{
\begin{array}{llll}
\para{i\partial_t+\Delta_\g+q(x)}v(t,x)=F(t,x)  & \textrm{in }\,\,(0,T)\times\M,\cr
\\
v(T,x)=0 & \textrm{in }\,\,\M,\cr
\\
v(t,x)=0 & \textrm{on } \,\, (0,T)\times\p\M.
\end{array}
\right.
\end{equation}
\end{defi}
One gets the following lemma.
\begin{Lemm}\label{L7.4}
Let $f\in L^2((0,T)\times\p\M)$. There exists a unique solution
\begin{equation}\label{7.30}
u\in\mathcal{C}(0,T;H^{-1}(\M))\cap H^{-1}(0,T;L^2(\M))
\end{equation}
defined by transposition, of the problem
\begin{equation}\label{7.31}
\left\{
\begin{array}{llll}
\para{i\partial_t+\Delta_\g+q(x)}u(t,x)=0  & \textrm{in }\,\,(0,T)\times\M,\cr
\\
u(0,x)=0 & \textrm{in }\,\,\M,\cr
\\
u(t,x)=f(t,x) & \textrm{on } \,\, (0,T)\times\p\M.
\end{array}
\right.
\end{equation}
Furthermore, there is a constant $C>0$ such that
\begin{equation}\label{7.32}
\norm{u}_{\mathcal{C}(0,T;H^{-1}(\M))}+\norm{u}_{H^{-1}(0,T;L^2(\M))}\leq C\norm{f}_{L^2((0,T)\times\p\M)}.
\end{equation}
\end{Lemm}
\begin{Demo}{}
Let $F\in\mathscr{H}=L^1(0,T;H^1_0(\M))$ or $\mathscr{H}=H_0^1(0,T;L^2(\M))$. Let $v\in\mathcal{C}(0,T;H^1_0(\M))$ solution of the backward
boundary value problem for the Schrödinger equation (\ref{7.33}). By Lemma \ref{L7.3} the mapping $F\mapsto\frac{\p v}{\p\nu}$ is linear and continuous
from $\mathscr{H}$ to $L^2((0,T)\times\M)$ and there exists $C>0$ such that
\begin{equation}\label{7.34}
\norm{v}_{\mathcal{C}(0,T;H^1_0(\M))}\leq C\norm{F}_{\mathscr{H}}
\end{equation}
and
\begin{equation}\label{7.35}
\norm{\p_\nu v}_{L^2((0,T)\times\p\M)}\leq C\norm{F}_{\mathscr{H}}.
\end{equation}
We define a linear functional $\ell$ on the linear space $\mathscr{H}$ as follows:
$$
\ell(F)=\seq{f,\p_\nu v}_0
$$
where $v$ solves (\ref{7.33}). By (\ref{7.35}), we obtain
$$
\abs{\ell(F)}\leq \norm{f}_{L^2((0,T)\times\p\M)}\norm{F}_{\mathscr{H}}.
$$
It is known that any linear bounded functional on the space $\mathscr{H}$ can be written as
$$
\ell(F)=\para{u,F}_{\mathscr{H}',\mathscr{H}}
$$
where $u$ is some element from the space $\mathscr{H}'$. Thus the system (\ref{7.31}) admits a solution $u\in\mathscr{H}'$ in the transposition sense,
which satisfies
$$
\norm{u}_{\mathscr{H}'}\leq C\norm{f}_{L^2((0,T)\times\p\M)}.
$$
This completes the proof of the Lemma.
\end{Demo}
In what follows, we will need the following estimate for non-homogenous elliptic boundary value problem.
\begin{Lemm}\label{L7.5}
Let $\psi\in H^{-1}(\M)$ and $\phi\in H^1(\p\M)$. Let $w\in H^1(\M)$ the unique solution of the following boundary value problem
\begin{equation}\label{7.36}
\left\{
\begin{array}{llll}
\Delta_\g w(x)=\psi(x)  & \textrm{in }\,\,\M,\cr
\\
w(x)=\phi & \textrm{on } \,\,\p\M,
\end{array}
\right.
\end{equation}
then the following estimate holds true
\begin{equation}\label{7.37}
\norm{w}_{H^1(\M)}\leq C\para{\norm{\psi}_{H^{-1}(\M)}+\norm{\phi}_{H^1(\p\M)}}.
\end{equation}
\end{Lemm}
\begin{Demo}{}
We decompose the solution $w$ of (\ref{7.36}) as $w=w_1+w_2$ with $w_1$ and $w_2$, respectively, solution of
\begin{equation}\label{7.38}
\left\{
\begin{array}{llll}
\Delta_\g w_1(x)=\psi(x)  & \textrm{in }\,\,\M,\cr
\\
w_1(x)=0 & \textrm{on } \,\,\p\M,
\end{array}
\right.
;\quad
\left\{
\begin{array}{llll}
\Delta_\g w_2(x)=0  & \textrm{in }\,\,\M,\cr
\\
w_2(x)=\phi & \textrm{on } \,\,\p\M,
\end{array}
\right.
\end{equation}
Since $-\Delta_\g$ is an isomorphism from $H^1_0(\M)$ to $H^{-1}(\M)$, we have
\begin{equation}\label{7.39}
\norm{w_1}_{H^1(\M)}\leq C\norm{\psi}_{H^{-1}(\M)}.
\end{equation}
Next, it is well known that (see \cite{[Necas]})
\begin{equation}\label{7.40}
\norm{w_2}_{L^2(\M)}\leq C\norm{\phi}_{L^2(\p\M)}.
\end{equation}
Now, we shall show that
\begin{equation}\label{7.41}
\norm{w_2}_{H^1(\M)}\leq C\norm{\phi}_{H^{1}(\p\M)}.
\end{equation}
Indeed, let $h\in H^1(\p\M)$ and $\theta\in\mathcal{C}_0^\infty(0,T)$, $\theta\ge 0$. Let $v$ solve the following initial boundary value problem
for the wave equation
\begin{equation}\label{7.42}
\left\{
\begin{array}{llll}
\para{\partial^2_t-\Delta_\g}v(t,x)=0  & \textrm{in }\,\,(0,T)\times\M,\cr
\\
v(0,x)=\p_tv(0,x)=0 & \textrm{in }\,\,\M,\cr
\\
v(t,x)=h(x)\theta(t) & \textrm{on } \,\, (0,T)\times\p\M.
\end{array}
\right.
\end{equation}
Then we have (see \cite{[KKL]})
$$
v\in\mathcal{C}(0,T;H^1(\M))\cap\mathcal{C}^1(0,T;L^2(\M)).
$$
Furthermore there exist $C>0$ such that
\begin{multline}\label{7.43}
\norm{v}_{\mathcal{C}(0,T;H^1(\M))}+\norm{v}_{\mathcal{C}^1(0,T;L^2(\M))}+\norm{\p_\nu v}_{L^2((0,T)\times\p\M)} \\ 
\leq C\norm{h}_{H^1(\p\M)}.
\end{multline}
Multiplying both sides of (\ref{7.42}) by $w_2$ and integrating over $(0,T)\times\M$, we get
\begin{align}\label{7.44}
0&=\int_0^T\!\!\!\int_\M\para{\p_t^2v-\Delta_\g v}w_2(x)\dv \, \dd t \cr &=\int_\M \p_tv(T,x) w_2(x)\dv
-\int_0^T\!\!\!\int_{\p\M}\p_\nu v \phi(x)\ds \, \dd t \cr &\quad+\int_0^T\!\!\!\theta(t) \, \dd t \int_{\p\M}h(x)\p_\nu w_2\ds.
\end{align}
Then, by (\ref{7.43}) and (\ref{7.40}), one gets
\begin{align}\label{7.45}
\abs{\int_{\p\M}h(x)\p_\nu w_2\ds}  &\leq C\big(\norm{\p_\nu v}_{L^2((0,T)\times\p\M)}\norm{\phi}_{L^2(\p\M)}\cr
&\quad +\norm{w_2}_{L^2(\M)}\norm{v}_{\mathcal{C}^1(0,T;L^2(\M))}\big)\cr
& \leq C\norm{\phi}_{L^2(\p\M)}\norm{h}_{H^1(\p\M)}.
\end{align}
This implies
\begin{equation}\label{7.46}
\norm{\p_\nu w_2}_{H^{-1}(\p\M)}\leq C\norm{\phi}_{L^2(\p\M)}.
\end{equation}
Furthermore, Green's formula yields
\begin{equation}\label{7.47}
\int_\M\abs{\nabla_\g w_2}^2\dv\leq \norm{\p_\nu w_2}_{H^{-1}(\p\M)}\norm{\phi}_{H^1(\p\M)}\leq C\norm{\phi}^2_{H^1(\p\M)}.
\end{equation}
>From (\ref{7.47}) and (\ref{7.40}), we get
\begin{equation}\label{7.48}
\norm{w_2}_{H^1(\M)}\leq C\norm{\phi}_{H^1(\p\M)}.
\end{equation}
Both (\ref{7.48}) and (\ref{7.39}) yield likewise
\begin{equation}\label{7.49}
\norm{w}_{H^1(\M)}\leq C\big(\norm{\psi}_{H^{-1}(\M)}+\norm{\phi}_{H^1(\p\M)}\big).
\end{equation}
This completes the proof of (\ref{7.37}).
\end{Demo}
\subsection{Proof of Theorem \ref{Th0}}
We proceed to prove Theorem \ref{Th0}. Let $f\in H^1((0,T)\times\p\M)$ and $u$ solve (\ref{1.2}). Put $u'=\p_tu$, then
\begin{equation}\label{7.50}
\left\{
\begin{array}{llll}
\para{i\partial_t+\Delta_\g+q(x)}u'(t,x)=0  & \textrm{in }\,\,(0,T)\times\M,\cr
\\
u'(0,x)=0 & \textrm{in }\,\,\M,\cr
\\
u'(t,x)=f'(t,x) & \textrm{on } \,\, (0,T)\times\p\M,
\end{array}
\right.
\end{equation}
Since $f'\in L^2((0,T)\times\p\M)$, by lemma \ref{L7.4}, we get
\begin{equation}\label{7.51}
u'\in\mathcal{C}(0,T;H^{-1}(\M))\cap H^{-1}(0,T;L^2(\M)).
\end{equation}
Furthermore there is a constant $C>0$ such that
\begin{equation}\label{7.52}
\norm{u'}_{\mathcal{C}(0,T;H^{-1}(\M))}+\norm{u'}_{H^{-1}(0,T;L^2(\M))}\leq C\norm{f}_{H^1((0,T)\times\p\M)}.
\end{equation}
Thus (\ref{7.51}) implies the following regularity for $u$
\begin{align*}
u &\in\mathcal{C}^1(0,T;H^{-1}(\M))\cap \mathcal{C}(0,T;L^2(\M)), \\
\Delta_\g u &\in\mathcal{C}(0,T;H^{-1}(\M))\cap H^{-1}(0,T;L^2(\M)).
\end{align*}
Since $f(t,\cdot)\in H^1(\p\M)$, by the elliptic regularity, we get
$$
u\in \mathcal{C}(0,T;H^1(\M))\cap \mathcal{C}^1(0,T;H^{-1}(\M)).
$$
Moreover there exists $C>0$ such that the following estimates hold true
\begin{align}\label{7.53}
&&\norm{u}_{\mathcal{C}^1(0,T;H^{-1}(\M))}\leq C\norm{f}_{H^1((0,T)\times\p\M)},\cr
&&\norm{\Delta_\g u}_{\mathcal{C}(0,T;H^{-1}(\M))}\leq C\norm{f}_{H^1((0,T)\times\p\M)}.
\end{align}
Using Lemma \ref{L7.5}, we find
\begin{align}\label{7.54}
&&\norm{u}_{\mathcal{C}^1(0,T;H^{-1}(\M))}\leq C\norm{f}_{H^1((0,T)\times\p\M)},\cr
&&\norm{u}_{\mathcal{C}(0,T;H^{1}(\M))}\leq C\norm{f}_{H^1((0,T)\times\p\M)}.
\end{align}
The proof of (\ref{1.8}) is as in Lemma \ref{L7.3}. If one multiplies (\ref{1.2}) by $\seq{N,\nabla_\g\overline{u}}_\g$, the arguments leading to
(\ref{7.20}) give now
\begin{align}\label{7.55}
0=&\int_0^T\!\!\!\int_\M i\p_tu\seq{N,\nabla_\g\overline{u}}_\g\dv \, \dd t+\int_0^T\!\!\!\int_\M\Delta_\g u
\seq{N,\nabla_\g\overline{u}}_\g\dv \, \dd t\cr
&+\int_0^T\!\!\!\int_\M q(x)u\seq{N,\nabla_\g\overline{u}}_\g\dv \, \dd t =I'_1+I'_2+I'_3,
\end{align}
with
\begin{equation}\label{7.56}
\abs{\textrm{Re}\, I'_1}\leq C_\varepsilon\norm{f}_{H^1((0,T)\times\p\M)}^2+\varepsilon\norm{\p_\nu u}^2_{L^2((0,T)\times\p\M)},
\end{equation}
where we have used (\ref{7.54}) instead of (\ref{7.13})-(\ref{7.7}). Furthermore, we derive from Green's formula
\begin{align}\label{7.57}
\textrm{Re}\,I'_2&=\frac{1}{2}\int_0^T\!\!\!\int_{\p\M}\abs{\p_\nu u}^2\ds \, \dd t  +\int_0^T\!\!\!\int_\M D_{\g}N(\nabla_\g u,
\nabla_\g\overline{u})\dv \, \dd t\cr
&\quad -\frac{1}{2}\int_0^T\!\!\!\int_\M\!\!\abs{\nabla_\g u}_\g^2\dive(N)\dv \, \dd t-\frac{1}{2}\int_0^T\!\!\!\int_{\p\M}\!\!\abs{\nabla_\tau f}^2\ds
\, \dd t.
\end{align}
This together with
\begin{equation}\label{7.58}
\abs{\textrm{Re}\, I'_3}\leq \norm{f}^2_{H^1((0,T)\times\p\M)}
\end{equation}
and (\ref{7.58}), (\ref{7.57}) and (\ref{7.56}) imply
\begin{equation}
\norm{\p_\nu u}_{L^2((0,T)\times\p\M)}\leq C\norm{f}_{H^1((0,T)\times\p\M)},
\end{equation}
where we have used (\ref{7.54}) again. The proof of Theorem \ref{Th0} is now complete.
\section{Geometrical optics solutions of the Schr\"odinger equation}
\setcounter{equation}{0}
We now proceed to the construction of geometrical optics solutions to the  Schr\"odinger equation. We extend the  manifold $(\M,\g)$ into a simple
manifold $\M_2 \Supset \M$ and consider a simple manifold $(\M_1,\g)$ such that $\M_2\Supset \M_1$. The potentials
$q_{1},q_{2}$ may also be extended to $\M_{2}$ and their $H^1(\M_{1})$ norms may be bounded by $M_{0}$. Since $q_{1}$ and $q_{2}$
coincide on the boundary, their extension outside $\M$ can be taken the same so that $q_{1}=q_{2}$ in $\M_{2} \setminus \M_{1}$.
\medskip

Our construction here is a modification of a similar result in \cite{[BellDSF]}, which dealt with the situation of the wave equation.
\medskip

We suppose, for a moment, that is able to find a function $\psi\in{\cal C}^2(\M)$ which satisfies the eikonal equation
\begin{equation}\label{4.1}
\abs{\nabla_\g\psi}^2_\g=\sum_{i,j=1}^n\g^{ij}\frac{\p\psi}{\p x_i}\frac{\p\psi}{\p
x_j}=1,\qquad \forall x\in \M_2
\end{equation}
and assume that there exist a function $a\in H^1(\R,H^2(\M))$ which solves the transport equation
\begin{equation}\label{4.2}
\frac{\p a}{\p t}+\sum_{j,k=1}^n \g^{jk}\frac{\p\psi}{\p x_j}\frac{\p
a}{\p x_k}+\frac{1}{2}(\Delta_\g\psi)a=0,\qquad \forall t\in\R,\, x\in\M
\end{equation}
with initial or final data
\begin{equation}\label{4.3}
a(t,x)=0,\quad \forall x\in \M,\quad \textrm{and}\,\,t\leq0,\,\,\textrm{or}\,\,t\geq T_0.
\end{equation}
We also introduce the norm $\norm{\cdot}_*$ given by
\begin{equation}\label{4.7}
\norm{a}_*=\norm{a}_{H^1(0,T_0;H^2(\M))}.
\end{equation}
\begin{Lemm}\label{L4.1}
Let $q\in L^\infty(\M)$. Then the following Schr\"odinger equation
\begin{align*}
(i\p_t+\Delta_\g +q(x))u&=0,\quad \textrm{in}\quad \M_T:=(0,T)\times\M,\\
u(\kappa,x)&=0,\quad\kappa=0,\textrm{ or } T
\end{align*}
has a solution of the form
\begin{equation}\label{4.4}
u(t,x)=a(2\lambda t,x)e^{i\lambda(\psi(x)-\lambda
t)}+v_\lambda(t,x),
\end{equation}
such that
\begin{equation}\label{4.5}
u\in {\cal C}^1(0,T;L^2(\M))\cap{\cal C}(0,T;H^2(\M))
\end{equation}
where $v_\lambda(t,x)$ satisfies
\begin{align*}
v_\lambda(t,x)&=0,\quad\forall (t,x)\in (0,T)\times\p\M , \\
v_\lambda(\kappa,x)&=0,\quad x\in \M,\quad \kappa=0 \textrm{ or } T.
\end{align*}
Furthermore, there exist $C>0$ such that, for any $\lambda>0$ the following estimates hold true.
\begin{equation}\label{4.6}
\norm{v_\lambda(t,\cdot)}_{H^k(\M)}\leq C\lambda^{k-1}\norm{a}_*,\qquad
k=0,1.
\end{equation}
The constant $C$ depends only on $T$ and $\M$ (that is $C$ does
not depend on $a$ and $\lambda$).
\end{Lemm}
\begin{Demo}{}
Let us consider
\begin{equation}\label{4.9}
k(t,x)=-\para{i\partial_t+\Delta_\g+q}\para{a(2\lambda
t,x)e^{i\lambda(\psi-\lambda t)}}.
\end{equation}
Let $v$ solve the following homogenous boundary value problem
\begin{equation}\label{4.8}
\left\{
\begin{array}{llll}
\para{i\partial_t+\Delta_\g+q}v(t,x)=k(t,x)
& \textrm{in }\,\, (0,T)\times\M,\cr
\\
v(\kappa,x)=0,& \textrm{in
}\,\,\M,\,\tau=0,\,\textrm{or}\,\,T\cr
\\
v(t,x)=0 & \textrm{on} \,\, (0,T)\times\p\M,
\end{array}
\right.
\end{equation}
To prove our Lemma it would be enough to show that $v$ satisfies the estimates (\ref{4.6}). We shall prove the estimate
for $\kappa=0$, and the $\kappa=T$ case may be handled in a similar
fashion. By a simple computation, we have
\begin{align}\label{4.10}
-k(t,x)&=e^{i\lambda(\psi(x)-\lambda t)}\para{\Delta_\g+q(x)}\para{a(2\lambda t,x)}\cr
&\quad +2i\lambda e^{i\lambda(\psi(x)-\lambda t)}\para{\p_ta+\sum_{j,k=1}^n \g^{jk}\frac{\p\psi}{\p x_j}
\frac{\p a}{\p x_k}+\frac{a}{2}\Delta_\g\psi}(2\lambda t,x)\cr
&\quad+\lambda^2 a(2\lambda t,x) e^{i\lambda(\psi(x)-\lambda t)}\para{1-\sum_{j,k=1}^n\g^{jk}\frac{\p\psi}{\p x_j}
\frac{\p\psi}{\p x_k}}.
\end{align}
Taking into account  (\ref{4.1}) and (\ref{4.2}), the right-hand side of (\ref{4.10}) becomes
\begin{align}\label{4.11}
k(t,x)&=-e^{i\lambda(\psi(x)-\lambda t)}\para{\Delta_\g+q}\para{a(2\lambda t,x)}
\cr &\equiv-e^{i\lambda(\psi(x)-\lambda t)}k_0(2\lambda t,x).
\end{align}
Since $k_0\in H^1_0(0,T;L^2(\M))$, by Lemma \ref{L7.2}, we find
\begin{equation}\label{4.12}
v_\lambda\in {\cal C}^1(0,T;L^2(\M))\cap{\cal C}(0,T;H^2(\M)\cap
H^1_0(\M)).
\end{equation}
Furthermore, there is a constant $C>0$, such that
\begin{align}\label{4.13}
\norm{v_\lambda(t,\cdot)}_{L^2(\M)}&\leq C\int_0^T\norm{k_0(2\lambda
t,\cdot)}_{L^2(\M)}\, \dd t \displaybreak[1] \\ \nonumber &\leq
\frac{C}{\lambda}\int_\R\norm{k_0(s,\cdot)}_{L^2(\M)} \, \dd s \displaybreak[1] \\ \nonumber  &\leq
\frac{C}{\lambda}\norm{a}_*.
\end{align}
Moreover, we have
\begin{multline}\label{4.14}
\norm{\nabla v_\lambda (t,\cdot)}_{L^2(\M)} \\ 
\leq
C\eta\int_0^T\para{\lambda^2\norm{k_0(2\lambda
t,\cdot)}_{L^2(\M)}+\lambda \norm{\p_t k_0(2\lambda
t,\cdot)}_{L^2(\M)}}\, \dd t\\ +\eta^{-1}\int_0^T\norm{
k_0(2\lambda t,\cdot)}_{L^2(\M)}\, \dd t.
\end{multline}
Finally, choosing $\eta=\lambda^{-1}$, we obtain
\begin{align}\label{4.15}
\norm{\nabla v_\lambda(t,\cdot)}_{L^2(\M)} &\leq
C\para{\int_\R\norm{k_0(s ,\cdot)}_{L^2(\M)} \, \dd s+\int_\R\norm{\p_t k_0(s ,\cdot)}_{L^2(\M)} \, \dd s} \nonumber \\ &\leq  C\norm{a}_*.
\end{align}
Combining (\ref{4.15}) and (\ref{4.13}), we immediately deduce the estimate (\ref{4.6}).
\end{Demo}
\medskip

We will now construct the phase function $\psi$ solution to the eikonal equation (\ref{4.1}) and the amplitude $a$ solution to the
transport equation (\ref{4.2}).
\medskip

Let $y\in \p \M_1$. Denote points in $\M_1$ by $(r,\theta)$
where $(r,\theta)$ are polar normal coordinates in $\M_1$ with center
$y$. That is $x=\exp_{y}(r\theta)$ where $r>0$ and
$$
\theta\in S_{y}\M_1=\set{\xi\in T_{y}\M_1,\,\,\abs{\xi}_\g=1}.
$$
In these coordinates (which depend on the choice of $y$) the
metric takes the form
$$
\widetilde{\g}(r,\theta)=\dd r^2+\g_0(r,\theta)
$$
where $\g_0(r,\theta)$ is a smooth positive definite metric.
For any function $u$ compactly supported in $\M$, we set for $r>0$ and $\theta\in S_y\M_1$
$$
\widetilde{u}(r,\theta)=u(\exp_{y}(r\theta))
$$
where we have extended $u$ by $0$ outside $\M$.
An explicit solution to the eikonal equation (\ref{4.1}) is the geodesic distance function to $y \in \p \M_1$
\begin{equation}\label{4.16}
\psi(x)=d_\g(x,y).
\end{equation}
By the simplicity assumption, since $y\in \M_2\backslash\overline{\M}$, we have $\psi\in\cal{C}^\infty(\M)$ and
\begin{equation}\label{4.17}
\widetilde{\psi}(r,\theta)=r=d_\g(x,y).
\end{equation}
The next step is to solve the transport equation (\ref{4.2}). Recall that
if $f(r)$ is any function of the geodesic distance $r$, then
\begin{equation}\label{4.18}
\Delta_{\widetilde{\g}}f(r)=f''(r)+\frac{\alpha^{-1}}{2}\frac{\p
\alpha}{\p r}f'(r).
\end{equation}
Here $\alpha=\alpha(r,\theta)$ denotes the square of the volume element in geodesic polar coordinates.
The transport equation (\ref{4.2}) becomes
\begin{equation}\label{4.19}
\frac{\p \widetilde{a}}{\p t}+\frac{\p \widetilde{\psi}}{\p
r}\frac{\p \widetilde{a}}{\p
r}+\frac{1}{4}\widetilde{a}\alpha^{-1}\frac{\p \alpha}{\p r}\frac{\p
\widetilde{\psi}}{\p r}=0.
\end{equation}
Thus $\widetilde{a}$ satisfy
\begin{equation}\label{4.20}
\frac{\p \widetilde{a}}{\p t}+\frac{\p \widetilde{a}}{\p
r}+\frac{1}{4}\widetilde{a}\alpha^{-1}\frac{\p \alpha}{\p r}=0.
\end{equation}
Let $\phi\in{\cal C}_0^\infty(\R)$ and $b\in H^2(\p_+S\M)$. Let us write $\widetilde{a}$ in the form
\begin{equation}\label{4.21}
\widetilde{a}(t,r,\theta)=\alpha^{-1/4}\phi(t-r)b(y,\theta).
\end{equation}
Direct computations yields
\begin{equation}\label{4.22}
\frac{\p \widetilde{a}}{\p
t}(t,r,\theta)=\alpha^{-1/4}\phi'(t-r)b(y,\theta).
\end{equation}
and
\begin{equation}\label{4.23}
\frac{\p \widetilde{a}}{\p
r}(t,r,\theta)=-\frac{1}{4}\alpha^{-5/4}\frac{\p\alpha}{\p
r}\phi(t-r)b(y,\theta)-\alpha^{-1/4}\phi'(t-r)b(y,\theta).
\end{equation}
Finally, (\ref{4.23}) and (\ref{4.22}) yield
\begin{equation}\label{4.24}
\frac{\p \widetilde{a}}{\p t}(t,r,\theta)+\frac{\p \widetilde{a}}{\p
r}(t,r,\theta)=-\frac{1}{4}\alpha^{-1}\widetilde{a}(t,r,\theta)\frac{\p\alpha}{\p
r}.
\end{equation}
Now if we assume that $\mathrm{supp}(\phi)\subset(0,\varepsilon_0)$, $\varepsilon_0>0$ small, then for any $x=\exp_y(r\theta)\in \M$, it is easy to
see that $\widetilde{a}(t,r,\theta)=0$ if $t\leq 0$ and $t\geq T_0$ for some $T_0>0$ sufficiently large.
\section{Stable determination of the electric potential}
\setcounter{equation}{0}
In this section, we complete the proof of Theorem \ref{Th1}. We are going to use the geometrical
optics solutions constructed in the previous section; this will provide information on the geodesic ray transform of the difference of electric potentials.
\subsection{Preliminary estimates}
The main purpose of this section is to present a preliminary estimate, which relates the
difference of the potentials to the Dirichlet-to-Neumann map.
As before, we let $q_1,\,q_2\in\mathscr{Q}(M_0)$ such that $q_1=q_2$ on the boundary $\p\M$. We set
$$
q(x)=(q_1-q_2)(x).
$$
Recall that we have extended $q_{1},q_{2}$ as $H^1(\M_{2})$ in such a way that $q=0$ on $\M_{2} \setminus \M$.
\begin{Lemm}\label{L5.1}
There exists $C>0$ such that for any $a_1$, $a_2\in H^1(\R, H^2(\M))$
satisfying the transport equation (\ref{4.2}) with (\ref{4.3}), the following estimate holds true:
\begin{multline}\label{5.1}
\abs{\int_{0}^T\!\!\!\!\int_{\M}q(x)a_1(2\lambda
t,x)\overline{a}_2(2\lambda t,x)\,\dv \, \dd t } \\ 
\leq C\para{\lambda^{-2}+\norm{\Lambda_{\g,\,q_1}-\Lambda_{\g,\,q_2}}}\norm{a_1}_*\norm{a_2}_*
\end{multline}
for any sufficiently large $\lambda > 0$.
\end{Lemm}
\begin{Demo}{} First, if $a_2$ satisfies (\ref{4.2}), (\ref{4.3}) and $\lambda$ is large enough, Lemma \ref{L4.1} guarantees
the existence of a geometrical optics solution $u_2$
\begin{equation}\label{5.2}
u_2(t,x)=a_2(2\lambda t,x)e^{i\lambda(\psi(x)-\lambda
t)}+v_{2,\lambda}(t,x),
\end{equation}
to the Schr\"odinger equation corresponding to the electric potential $q_2$,
$$
\para{i\p_t+\Delta_{\g}+q_2(x)}u(t,x)=0\quad \textrm{in}\,(0,T)\times\M, \quad u(0,\cdot)=0 \quad \textrm{in}\, \M
$$
where
$v_{2,\lambda}$ satisfies
\begin{gather}\label{5.3}
\lambda\norm{v_{2,\lambda}(t,\cdot)}_{L^2(\M)}+\norm{\nabla
v_{2,\lambda}(t,\cdot)}_{L^2(\M)}\leq C\norm{a_2}_*
\\ \nonumber
v_{2,\lambda}(t,x)=0,\quad\forall (t,x)\in\,(0,T)\times\p\M.
\end{gather}
Moreover
$$
u_2\in  {\cal C}^1(0,T;L^2(\M))\cap {\cal C}(0,T;H^2(\M)).
$$
Let us denote by $f_\lambda$ the function
$$
f_\lambda(t,x)=a_2(2\lambda t,x)e^{i\lambda(\psi(x)-\lambda
t)},\quad  t\in(0,T),\,\,x\in \p \M.
$$
Let us consider $v$ the solution of the following non-homogenous boundary value problem
\begin{equation}\label{5.4}
\left\{\begin{array}{lll}
\para{i\p_t+\Delta_{\g}+q_1} v=0, & (t,x)\in (0,T)\times\M,\cr
\\
v(0,x)=0, & x\in \M,\cr
\\
v(t,x)=u_2(t,x):=f_{\lambda}(t,x), & (t,x)\in (0,T)\times\p\M.
\end{array}
\right.
\end{equation}
Denote $ w=v-u_2$. Therefore, $w$ solves the following homogenous boundary value problem for the Schr\"odinger equation
$$
\left\{\begin{array}{lll}
\para{i\p_t+\Delta_{\g}+q_1(x)}w(t,x)=q(x)u_2(t,x) & (t,x)\in (0,T)\times\M,\cr
\\
w(0,x)=0, & x\in \M,\cr
\\
w(t,x)=0, & (t,x)\in (0,T)\times\p\M.
\end{array}
\right.
$$
Using the fact that  $q(x)u_2\in W^{1,1}(0,T;L^2(\M))$ with $u(0,\cdot)\equiv 0$, by Lemma \ref{L7.2}, we deduce that
$$
w\in {\cal C}^1(0,T;L^2(\M))\cap {\cal C}(0,T;H^2(\M)\cap
H^1_0(\M)).
$$
Therefore, we have constructed a special solution
  $$ u_1\in {\cal C}^1(0,T;L^2(\M))\cap {\cal C}(0,T;H^2(\M)) $$
to the backward Schr\"odinger equation
\begin{align*}
\para{i\partial_t+\Delta_{\g}+q_1(x)}u_1(t,x)&=0,  \quad (t,x) \in (0,T)\times\M, \\
u_1(T,x)&=0,  \quad x \in \M,
\end{align*}
having the special form
\begin{equation}\label{5.5}
u_1(t,x)=a_1(2\lambda t,x)e^{i\lambda(\psi(x)-\lambda
t)}+v_{1,\lambda}(t,x),
\end{equation}
which corresponds to the electric potential $q_1$, where
$v_{1,\lambda}$ satisfies
\begin{equation}\label{5.6}
\lambda\norm{v_{1,\lambda}(t,\cdot)}_{L^2(\M)}+\norm{\nabla
v_{1,\lambda}(t,\cdot)}_{L^2(\M)}\leq C\norm{a_1}_*.
\end{equation}
Integrating by parts and using Green's formula (\ref{2.4}), we find
\begin{align}\label{5.7}
\int_0^T\!\!\!\int_\M\para{i\p_t+\Delta_\g+q_1}w\overline{u}_1\dv \, \dd t
&= \int_0^T\!\!\!\int_\M q u_2\overline{u}_1\dv \, \dd t\cr
&=-\int_0^T\!\!\!\int_{\p \M}\p_\nu w\overline{u}_1\ds \, \dd t.
\end{align}
Taking (\ref{5.7}), (\ref{5.4}) into account, we deduce
\begin{multline}\label{5.8}
\int_0^T\!\!\!\int_\M q(x)u_2(t,x)\overline{u}_1(t,x)\dv\,\dd t \\
=-\int_0^T\!\!\!\int_{\p \M}\para{\Lambda_{\g,\,q_1}-\Lambda_{\g,\,q_2}}(f_{\lambda})(t,x)\overline{g}_\lambda(t,x) \ds \, \dd t
\end{multline}
where $g_\lambda$ is given by
$$
g_\lambda(t,x)=a_1(2\lambda t,x)e^{i\lambda(\psi(x)-\lambda
t)},\quad (t,x)\in (0,T)\times\p\M.
$$
It follows from (\ref{5.8}), (\ref{5.5}) and (\ref{5.2}) that
\begin{align}\label{5.9}
\int_0^T\!\!\!\int_\M q(x) (a_2 \overline{a}_1)(2\lambda t,x)\dv\,\dd t &=
-\int_0^T\!\!\!\int_{\p \M} \overline{g}_\lambda \para{\Lambda_{\g,\,q_1}-\Lambda_{\g,\,q_2}}f_{\lambda} \ds \, \dd t \nonumber \\
&\quad-\int_0^T\!\!\!\int_\M q e^{i\lambda(\psi-\lambda t)}a_2(2\lambda t,x)\overline{v}_{1,\lambda}\dv \, \dd t \nonumber \\
&\quad-\int_0^T\!\!\!\int_\M q v_{2,\lambda}e^{-i\lambda(\psi-\lambda t)}\overline{a}_1(2\lambda t,x)\dv\, \dd t \nonumber \\
&\quad-\int_0^T\!\!\!\int_\M q v_{2,\lambda}\overline{v}_{1,\lambda}\dv \, \dd t.
\end{align}
In view of (\ref{5.6}) and (\ref{5.3}), we have
\begin{align}\label{5.10}
\bigg|\int_0^T\!\!\!\int_\M q e^{i\lambda(\psi-\lambda t)}&a_2(2\lambda t,x)\overline{v}_{1,\lambda}\dv\,\dd t\bigg|  \nonumber \\ &\leq
C\int_0^T\norm{a_2(2\lambda t,\cdot)}_{L^2(\M)}\norm{v_{1,\lambda}(t,\cdot)}_{L^2(\M)}\, \dd t \nonumber
\\ & \leq C\lambda^{-2}\norm{a_2}_*\norm{a_1}_*.
\end{align}
Similarly,  we deduce
\begin{equation}\label{5.10'}
\abs{\int_0^T\!\!\!\int_\M q e^{-i\lambda(\psi-\lambda t)}
\overline{a}_1(2\lambda t,x)v_{2,\lambda}(t,x)\dv\,\dd t}\leq C\lambda^{-2}\norm{a_1}_*\norm{a_2}_*.
\end{equation}
Moreover we have
\begin{equation}\label{5.10"}
\abs{\int_0^T\!\!\!\int_\M q(x)
v_{2,\lambda}(t,x)\overline{v}_{1,\lambda}(t,x)\dv\,\dd t}\leq C\lambda^{-2}\norm{a_1}_*\norm{a_2}_*.
\end{equation}
On the other hand, by the trace theorem, we find
\begin{align}\label{5.11}
\bigg|\int_0^T\!\!\!\int_{\p \M}&\para{\Lambda_{\g,\,q_1}-\Lambda_{\g,\,q_2}}(f_{\lambda}) \overline{g}_\lambda \ds \, \dd t \bigg| \cr
&\leq \norm{\Lambda_{\g,\,q_1}-\Lambda_{\g,\,q_2}} \norm{f_\lambda}_{H^1((0,T)\times\p\M)}\norm{g_\lambda}_{L^2((0,T)\times\p\M)}\cr
&\leq C\frac{\lambda^{1/2}}{\lambda^{1/2}}\norm{a_1}_*\norm{a_2}_*\norm{\Lambda_{\g,\,q_1}-\Lambda_{\g,\,q_2}}.
\end{align}
The estimate (\ref{5.1}) follows easily from (\ref{5.9}), (\ref{5.10}), (\ref{5.10'}), (\ref{5.10"}) and (\ref{5.11}).
This completes the proof of the Lemma.
\end{Demo}
\begin{Lemm}\label{L5.2} Let $M_0>0$. There exists $C>0$ such that for any $b\in H^2(\p_+S\M_{1})$, the following estimate
\begin{multline}
\abs{\int_{S_{y}\M_1}\!\int^{\tau_+(y,\theta)}_0\widetilde{q}(s,\theta)b(y,\theta) \mu(y,\theta) \, \dd s \, \dd \omega_y(\theta)} \\ \leq
C\norm{\Lambda_{\g,\,q_1}-\Lambda_{\g,\,q_2}}^{1/2} \norm{b(y,\cdot)}_{H^2(S_y^+\M_{1})}.
\end{multline}
holds for any $y\in\p\M_1$.
\end{Lemm}
We use the notation
    $$ S_y^+\M_{1} = \big\{\theta \in S_{y}\M_{1} : \langle \nu,\theta \rangle_{\g}<0 \big\}. $$
\begin{Demo}{}
Following (\ref{4.21}), we take two solutions to (\ref{4.2}) and (\ref{4.3}) of the form
\begin{align*}
\widetilde{a}_1(t,r,\theta)&=\alpha^{-1/4}\phi(t-r)b(y,\theta), \\
\widetilde{a}_2(t,r,\theta)1&=\alpha^{-1/4}\phi(t-r)\mu(y,\theta).
\end{align*}
We recall that $\mu(y,\theta)=\langle \nu(y),\theta \rangle$ is the density of the $L^2$ space where the image of the geodesic ray transform lies.
Now we change variable in (\ref{5.1}), $x=\exp_{y}(r\theta)$, $r>0$ and
$\theta\in S_{y}\M_1$, we have
\begin{align*}
\int_0^T\!\!\int_\M &q a_1(2\lambda t,x)a_2(2\lambda t,x) \dv \, \dd t \cr
&=\int_0^T\!\!\int_{S_{y}\M_1}\!\int_0^{\tau_+(y,\theta)}\widetilde{q}(r,\theta)\widetilde{a}_1(2\lambda t,r,\theta)\widetilde{a}_2(2\lambda t,r,\theta)
\alpha^{1/2} \, \dd r \, \dd\omega_y(\theta) \, \dd t\cr
&=\int_0^T\!\!\int_{S_{y}\M_1}\!\int_0^{\tau_+(y,\theta)}\widetilde{q}(r,\theta)\phi^2(2\lambda t-r)b(y,\theta) \mu(y,\theta) \, \dd r \,
\dd\omega_y(\theta)
\, \dd t\cr
&=\frac{1}{2\lambda}\int_0^{2\lambda T}\!\!\!\int_{S_{y}\M_1}\!\int_0^{\tau_+(y,\theta)}\widetilde{q}(r,\theta)\phi^2( t-r)b(y,\theta)
\mu(y,\theta) \, \dd r \, \dd\omega_y(\theta) \, \dd t.
\end{align*}
By virtue of Lemma \ref{L5.1}, we conclude that
\begin{multline}\label{5.13}
\abs{\int_0^{\infty}\!\!\!\int_{S_{y}\M_1}\!\int_0^{\tau_+(y,\theta)}\widetilde{q}(r,\theta)\phi^2( t-r)b(y,\theta) \mu(y,\theta) \, \dd r
\, \dd\omega_y(\theta) \, \dd t} \\
\leq C\para{\lambda^{-1}+\lambda\norm{\Lambda_{\g,\,q_1}-\Lambda_{\g,\,q_2}}}\norm{\phi}^2_{H^3(\R)}\norm{b(y,\cdot)}_{H^2(S^+_y\M_{1})}.
\end{multline}
By the support properties of the function $\phi$, we get that the left-hand side term in the previous inequality reads
\begin{multline*}
    \int_0^{\infty}\!\!\!\int_{S_{y}\M_1}\!\int_0^{\tau_+(y,\theta)}\widetilde{q}(r,\theta)\phi^2( t-r)b(y,\theta) \mu(y,\theta) \, \dd r
    \, \dd\omega_y(\theta) \, \dd t \\
    = \bigg(\int_{-\infty}^{\infty} \phi(t) \dd t \bigg) \!\!\!\int_{S_{y}\M_1}\!\int_0^{\tau_+(y,\theta)}\widetilde{q}(r,\theta)b(y,\theta)
    \mu(y,\theta) \, \dd r \, \dd\omega_y(\theta).
\end{multline*}
Finally, minimizing in $\lambda$ in the right hand-side of \eqref{5.13} we obtain
\begin{multline*}
\abs{\int_{S_{y}\M_1}\!\int_0^{\tau_+(y,\theta)}\widetilde{q}(s,\theta)b(y,\theta) \mu(y,\theta) \, \dd s \, \dd\omega_y(\theta)} \\ 
\leq
C\norm{\Lambda_{\g,\,q_1}-\Lambda_{\g,\,q_2}}^{1/2} \norm{b(y,\cdot)}_{H^2(S^+_y\M_{1})}.
\end{multline*}
This completes the proof of the lemma.
\end{Demo}
\subsection{End of the proof of the stability estimate}
Let us now complete the proof of the stability estimate in Theorem \ref{Th1}. Using Lemma \ref{L5.2}, for any $y\in\p \M_1$ and $b\in H^2(\p_+ S\M)$
we have
\begin{multline*}
\abs{\int_{S_{y}\M_1}\I(q)(y,\theta)b(y,\theta) \mu(y,\theta) \, \dd \omega_y(\theta)} \\
\leq C\norm{\Lambda_{\g,\,q_1}-\Lambda_{\g,\,q_2}}^{1/2}\norm{b(y,\cdot)}_{H^2(S_y^+\M_{1})}.
\end{multline*}
Integrating with respect to $y\in \p \M_1$ we obtain
\begin{multline}\label{5.18}
\abs{\int_{\p_+S\M_1}\I(q)(y,\theta)b(y,\theta)\seq{\theta,\nu(y)}\dss (y,\theta)} \\ 
\leq C\norm{\Lambda_{\g,\,q_1}-\Lambda_{\g,\,q_2}}^{1/2}
\norm{b}_{H^2(\p_+S\M_{1})}.
\end{multline}
Now we choose
$$
b(y,\theta)=\I\para{\I^*\I(q)}(y,\theta).
$$
Taking into account (\ref{3.8}) and (\ref{3.4}), we obtain
$$
\norm{\I^*\I(q)}^2_{L^2(\M_1)}\leq C\norm{\Lambda_{\g,\,q_1}-\Lambda_{\g,\,q_2}}^{1/2} \norm{q}_{H^1(\M)}.
$$
By interpolation, it follows that
\begin{align}\label{5.19}
\norm{\I^*\I(q)}^2_{H^1(\M_1)}&\leq C \norm{\I^*\I(q)}_{L^2(\M_1)}\norm{\I^*\I(q)}_{H^2(\M_1)}\cr
&\leq C \norm{\I^*\I(q)}_{L^2(\M_1)}\norm{q}_{H^1(\M)}\cr
&\leq C \norm{\I^*\I(q)}_{L^2(\M_1)}\cr
&\leq C\norm{\Lambda_{\g,\,q_1}-\Lambda_{\g,\,q_2}}^{1/4}.
\end{align}
Using (\ref{3.7}), we deduce that
$$
\norm{q}_{L^2(\M)}^2\leq C\norm{\Lambda_{\g,\,q_1}-\Lambda_{\g,\,q_2}}^{1/4}.
$$
This completes the proof of Theorem \ref{Th1}.
\section{Stable determination of the conformal factor}
\setcounter{equation}{0}
This section is devoted to the proof of the stability estimate for the conformal factor. We use the following notations; let
$c\in\mathscr{C}(M_0,k,\varepsilon)$, we denote
\begin{gather}\label{6.1}
\varrho_0(x)=1-c(x),\quad\varrho_1(x)=c^{n/2}(x)-1,\quad\varrho_2(x)=c^{n/2-1}(x)-1,\cr
\varrho(x)=\varrho_2(x)-\varrho_1(x)=c^{n/2-1}(x)(1-c(x)).
\end{gather}
Then the following holds
\begin{gather}\label{6.2}
\norm{\varrho_j}_{\mathcal{C}^1(\M)}\leq C\norm{\varrho_0}_{\mathcal{C}^1(\M)},\quad j=1,2 \cr
C^{-1}\norm{\varrho_0}_{L^2(\M)}\leq\norm{\varrho}_{L^2(\M)}\leq C\norm{\varrho_0}_{L^2(\M)}.
\end{gather}
As in the case of potentials, we extend the  manifold $(\M,\g)$ into a simple manifold
$\M_2 \Supset \M$ so that $\M_{2} \Supset \M_{1} \Supset \M$ with $(\M_1,\g)$ simple. We extend the conformal factor $c$
by $1$ outside the manifold $\M$; its $\mathcal{C}^k(\M_{1})$ norms  may also be bounded by $M_{0}$.
The first step in our analysis is the following lemma.
\begin{Lemm}\label{L6.1}
Let $c\in \mathcal{C}^\infty(\M)$ such that $c=1$ on the boundary $\p\M$. Let $u_1$, $u_2$ solve the following boundary problems in $(0,T)\times\M$ with
some $T>0$
\begin{align}\label{6.4}
&\left\{
\begin{array}{llll}
(i\partial_t+\Delta_\g)u_1=0,  & \textrm{in }\; (0,T)\times \M\cr
\\
u_1(0,\cdot )=0, & \textrm{in }\; \M\cr
\\
u_1=f_1, & \textrm{on} \,\,(0,T)\times\p \M
\end{array}
\right.
\\
&\left\{
\begin{array}{llll}
(i\partial_t+\Delta_{\cg})u_2=0,  & \textrm{in }\; (0,T)\times \M\cr
\\
u_2(0,\cdot )=0, & \textrm{in }\; \M\cr
\\
u_2=f_2, & \textrm{on} \,\,(0,T)\times\p \M
\end{array}
\right.
\end{align}
Then the following identity
\begin{multline}\label{6.5}
\int_{0}^T\!\!\!\int_{\p\M}\para{\Lambda_\g-\Lambda_{\cg}}f_1\,\overline{f}_2\,\ds \, \dd t
=i\int_0^T\!\!\!\int_\M\varrho_1(x)u_1\p_t\overline{u}_2\, \dv \dd t \cr +\int_0^T\!\!\!\int_\M \varrho_2(x)\seq{\nabla_\g u_1(t,x),\nabla_\g
\overline{u}_2(t,x)}_\g\dv \,  \dd t
\end{multline}
holds true for any $f_j\in H^1((0,T)\times\p\M)$, $j=1,2$.
\end{Lemm}
\begin{Demo}{}
We multiply both hand sides of the first equation (\ref{6.4}) by $\overline{u}_2$, integrate by parts in time and use Green's formula (\ref{2.4}) to get
\begin{align*}
0&=\int_0^T\!\!\!\int_\M\para{i\p_tu_1+\Delta_\g u_1}\overline{u}_2\dv \, \dd t  \displaybreak[1] \\
&=-i\int_0^T\!\!\!\int_\M u_1\p_t\overline{u_2}\dvc \dd t+i\int_0^T\!\!\!\int_\M \varrho_1 u_1\p_t\overline{u}_2\dv \, \dd t \displaybreak[1] \\
&\quad+\int_0^T\!\!\!\int_{\p\M}\p_\nu u_1\overline{f_2} \ds \dd t -\int_0^T\!\!\!\int_\M\sum_{j,k=1}^n c\g^{jk}\para{\frac{\p u_1}{\p x_j}
\frac{\p \overline{u_2}}{\p x_k}}\dvc \, \dd t \displaybreak[1] \\
&\quad+\int_0^T\!\!\!\int_\M\varrho_2\para{\sum_{j,k=1}^n c\g^{jk}\frac{\p u_1}{\p x_j}\frac{\p \overline{u_2}}{\p x_k}}\dv \, \dd t
\end{align*}
and after using a second time Green's formula, we end up with
\begin{align*}
0&=i\int_0^T\!\!\!\int_\M \varrho_1 u_1\p_t\overline{u_2} \dv \, \dd t+\int_0^T\!\!\!\int_M\varrho_2 \para{\sum_{j,k=1}^n c\g^{jk}\frac{\p u_1}{\p x_j}
\frac{\p \overline{u_2}}{\p x_k}}\dv \, \dd t \displaybreak[1] \\
&\quad+\int_0^T\!\!\!\int_\M u_1\para{-i\p_t\overline{u_{2}}+\Delta_{\cg}\overline{u_2}}\dvc \, \dd t-\int_0^T\!\!\!\int_{\p\M}\p_\nu \overline{u_2}
f_1\dsc \, \dd t\\
&\quad+\int_0^T\!\!\!\int_{\p\M}\p_\nu u_1\overline{f_2} \,\ds \, \dd t.
\end{align*}
Taking into account the fact that $c=1$ on $\p\M$, the fact that $\para{-i\p_t\overline{u}_2+\Delta_{\cg}\overline{u}_2}=0$ in $(0,T)\times\M$, and the
fact that the Dirichlet-to-Neumann map is selfadjoint, it follows that
\begin{multline}\label{6.7}
\int_{0}^T\!\!\!\int_{\p\M}\para{\Lambda_\g-\Lambda_{\cg}}f_1\,\overline{f_2} \ds \, \dd t =i\int_0^T\!\!\!\int_\M\varrho_1(x)u_1\p_t\overline{u_2}
\,\dv \, \dd t \\
+\int_0^T\!\!\!\int_\M \varrho_2(x)\para{\sum_{j,k=1}^n \g^{jk}\frac{\p u_1}{\p x_j}\frac{\p \overline{u_2}}{\p x_k}}\dv \, \dd t
\end{multline}
This completes the proof of the Lemma.
\end{Demo}
\subsection{Modified geometrical optics solutions}
Let $\psi_1$, $\psi_2$ be two phase functions solving the eikonal equation with respect to the metrics $\g$ and $\cg$.
\begin{equation}\label{6.8}
\begin{split}
\abs{\nabla_\g\psi_1}^2_\g&=\sum_{j,k=1}^n\g^{jk}\frac{\p\psi_1}{\p x_j}\frac{\p\psi_1}{\p x_k}=1,\\
\abs{\nabla_{\cg}\psi_2}^2_{\cg}&=\sum_{j,k=1}^n \cg^{jk}\frac{\p\psi_2}{\p x_j}\frac{\p\psi_2}{\p x_k}=1,
\end{split}
\quad\textrm{ on } \M.
\end{equation}
Let $a_2$ solve the transport equation in $\R\times\M$ with respect the metric $\g$ (as given in section 4)
\begin{equation}\label{6.9}
\frac{\p a_2}{\p t}+\sum_{j,k=1}^n\g^{jk}\frac{\p\psi_1}{\p x_j}\frac{\p a_2}{\p x_k}+\frac{a_2}{2}\Delta_\g\psi_1=0.
\end{equation}
Let $a_3$ solve the following transport equation in $\R\times\M$ with respect to the metric~$\cg$
\begin{align}\label{6.10}
\nonumber
\frac{\p a_3}{\p t}+\sum_{j,k=1}^n(\cg)^{jk}\frac{\p\psi_2}{\p x_j}\frac{\p a_3}{\p x_k}+\frac{a_3}{2}\Delta_{\cg}\psi_2&=-
\frac{1}{2i}a_2(t,x)(1-c^{-1})e^{i\lambda(\psi_1-\psi_2)} \\ &\equiv a_2(t,x)\varphi_0(x,\lambda)
\end{align}
which satisfies the bound
\begin{equation}\label{6.11}
\norm{a_3}_*\leq C\lambda \norm{\varrho_{0}}_{\mathcal{C}^1(\M)} \norm{a_2}_*.
\end{equation}
Let us now explain how to construct a solution $a_3$ satisfying (\ref{6.10}) and (\ref{6.11}). To solve the transport equation (\ref{6.10})
and (\ref{6.11}) it is enough to take, in the geodesic polar coordinates $(r,\theta)$ (with respect to the metric $\cg$)
\begin{equation}\label{6.12}
\widetilde{a}_3(t,r,\theta;\lambda)=\alpha_{\cg}^{-1/4}(r,\theta)\int_0^r\alpha_{\cg}^{1/4}(s,\theta)\widetilde{a}_2(s-r+t,s,\theta)
\widetilde{\varphi}_0(s,\theta,\lambda) \, \dd s,
\end{equation}
where $\alpha_{\cg}(r,\theta)$ denotes the square of the volume element in geodesic polar coordinates with respect to the metric $\cg$.
Using that $\norm{\varphi_0(\cdot,\lambda)}_{\mathcal{C}^1(\M)}\leq C\lambda  \norm{\varrho_{0}}_{\mathcal{C}^1(\M)} $ and (\ref{6.12})
we obtain (\ref{6.11}).
\begin{Lemm}\label{L6.2}
Let $c\in \mathscr{C}(M_0,k,\varepsilon)$ be such that $c=1$ near the boundary $\p\M$. Then the equation
\begin{equation}\label{6.13}
\para{i\p_t+\Delta_{\cg}}u=0,\quad\textrm{in}\quad (0,T)\times\M,\quad u(0,x)=0
\end{equation}
has a solution of the form
\begin{equation}\label{6.14}
u_2(t,x)=\frac{1}{\lambda}a_2(2\lambda t,x)e^{i\lambda(\psi_1(x)-\lambda t)}+a_3(2\lambda t,x;\lambda)e^{i\lambda(\psi_2(x)-\lambda t)}+
v_{2,\lambda}(t,x)
\end{equation}
which satisfies
\begin{multline}\label{6.15}
\lambda\norm{v_{2,\lambda}(t,\cdot)}_{L^2(\M)}+\norm{\nabla v_{2,\lambda}(t,\cdot)}_{L^2(\M)}+\lambda^{-1}\norm{\p_tv_{2,\lambda}(t,\cdot)}_{L^2(\M)}
\\ \leq C\para{\lambda^2\norm{\varrho_{0}}_{\mathcal{C}^1(M)}+\lambda^{-1}}\norm{a_2}_*
\end{multline}
where the constant $C$ depends only on $T$ and $\M$ (that is $C$ does not depend on $a$, $\lambda$ and $\varepsilon$).
\end{Lemm}
\begin{Demo}{}
We set
\begin{equation*}
k(t,x)=-\para{i\partial_t+\Delta_{\cg}}\para{\frac{1}{\lambda}a(2\lambda t,x)e^{i\lambda(\psi_1-\lambda t)}
-a_3(2\lambda t,x,\lambda)e^{i\lambda(\psi_2-\lambda t)}}.
\end{equation*}
To prove our Lemma it is enough to show that if $v$ solves
\begin{equation}\label{6.16}
\para{i\partial_t+\Delta_{\cg}}v=k
\end{equation}
with initial and boundary conditions
\begin{equation}\label{6.17}
v(0,x)=0,\quad\textrm{in}\quad \M,\quad \textrm{and}\quad v(t,x)=0\quad\textrm{on}\quad(0,T)\times\p\M
\end{equation}
then the estimate (\ref{6.15}) holds. But we have
\begin{align}\label{6.19}
-k(t,x)&=\frac{1}{\lambda}e^{i\lambda(\psi_1-\lambda t)}\Delta_{\cg}a_2(2\lambda t,x) \displaybreak[1]  \\ \nonumber
&\quad+2i e^{i\lambda(\psi_1-\lambda t)}\para{\p_ta_2+\sum_{j,k=1}^n \cg^{jk}\frac{\p\psi_1}{\p x_j}\frac{\p a_2}{\p x_k}+
\frac{a_2}{2}\Delta_{\cg}\psi_1}(2\lambda t,x) \displaybreak[1] \\  \nonumber
&\quad+\lambda e^{i\lambda(\psi_1-\lambda t)}a_2(2\lambda t,x)\para{1-c^{-1}\sum_{j,k=1}^n\g^{jk}\frac{\p\psi_1}{\p x_j}\frac{\p\psi_1}{\p x_k}}\cr
&\quad +e^{i\lambda(\psi_2-\lambda t)}\para{\Delta_{\cg}}\para{a_3(2\lambda t,x)} \displaybreak[1] \\  \nonumber
&\quad +2i\lambda e^{i\lambda(\psi_2-\lambda t)}\para{\p_ta_3+\sum_{j,k=1}^n\cg^{jk}\frac{\p\psi_2}{\p x_j}\frac{\p a_3}{\p x_k}
+\frac{a_3}{2}\Delta_{\cg}\psi_2}(2\lambda t,x) \displaybreak[1] \\  \nonumber
&\quad+\lambda^2e^{i\lambda(\psi_2-\lambda t)}a_3(2\lambda t,x) \para{1-\sum_{j,k=1}^n(\cg)^{jk}\frac{\p\psi_2}{\p x_j}\frac{\p\psi_2}{\p x_k}}.
\end{align}
Taking into account (\ref{6.8}) and (\ref{6.9}), the right-hand side of (\ref{6.19}) becomes
\begin{align}\label{6.20}
-k(t,x)&=\frac{1}{\lambda}e^{i\lambda(\psi_1-\lambda
t)}\Delta_{\cg} a_2(2\lambda t,x) \displaybreak[1] \\ \nonumber
&\quad
\begin{aligned}
    +2ie^{i\lambda(\psi_1-\lambda t)}&\bigg((c^{-1}-1)\seq{\nabla_\g\psi_1,\nabla_\g a_2(2\lambda t,x)}_\g \\
    &\quad+\frac{1}{2}a_2(2\lambda t,x)\para{\Delta_{\cg}\psi_1-\Delta_\g\psi_1}\bigg)
\end{aligned} \displaybreak[1] \\ \nonumber
&\quad
\begin{aligned}
   +2i\lambda e^{i\lambda(\psi_2-\lambda t)}&\bigg(\p_ta_3+\sum_{j,k=1}^n\cg^{jk}\frac{\p\psi_2}{\p x_j}\frac{\p a_3}{\p x_k}
   +\frac{a_3}{2}\Delta_{\cg}\psi_2 \\
   &\quad+\frac{a_2}{2i}e^{i\lambda(\psi_1-\psi_2)}(1-c^{-1})\bigg)(2\lambda t,x)
\end{aligned} \displaybreak[1] \\ \nonumber
&\quad+e^{i\lambda(\psi_2-\lambda t)}\Delta_{\cg}a_3(2\lambda t,x).
\end{align}
By (\ref{6.10}) we get
\begin{align*}
-k(t,x)&=\frac{1}{\lambda}e^{i\lambda(\psi_1-\lambda t)}\Delta_{\cg} a_2(2\lambda t,x)\cr
&\quad
\begin{aligned}
     +2ie^{i\lambda(\psi_1-\lambda t)}&\bigg((c^{-1}-1)\seq{\nabla_\g\psi_1,\nabla_\g a_2(2\lambda t,x)}_\g \\
     &\quad+\frac{1}{2}a_2(2\lambda t,x)\para{\Delta_{\cg}\psi_1-\Delta_\g\psi_1}\bigg)
\end{aligned}
\cr
&\quad+e^{i\lambda(\psi_2-\lambda t)}\Delta_{\cg}a_3(2\lambda t,x) \cr
&\equiv \frac{1}{\lambda}e^{i\lambda(\psi_1-\lambda t)}k_0(2\lambda t,x)+e^{i\lambda(\psi_1-\lambda t)}k_1(2\lambda t,x)
+e^{i\lambda(\psi_2-\lambda t)}k_2(2\lambda t,x).
\end{align*}
Since $k_j\in H^1_0(0,T;L^2(\M))$, by Lemma \ref{L7.2}, we deduce that
\begin{equation}\label{6.22}
v_\lambda\in {\cal C}^1(0,T;L^2(\M))\cap{\cal C}(0,T;H^2(\M)\cap
H^1_0(\M))
\end{equation}
and
\begin{align*}
\|&v_\lambda (t,\cdot)\|_{L^2(\M)} \\
&\leq C\! \int_0^T\!\!\para{\frac{1}{\lambda}\norm{k_0(2\lambda t,\cdot)}_{L^2(\M)}+\norm{k_1(2\lambda t,\cdot)}_{L^2(\M)}
+\norm{k_2(2\lambda t,\cdot)}_{L^2(\M)} } \, \dd t \displaybreak[1] \\  \nonumber
&\leq \frac{C}{\lambda}\int_\R\para{\frac{1}{\lambda}\norm{k_0(s,\cdot)}_{L^2(\M)}+\norm{k_1(s,\cdot)}_{L^2(\M)}+\norm{k_2(s,\cdot)}_{L^2(\M)}} \, \dd s\\ \nonumber
&\leq \frac{C}{\lambda}\para{\frac{1}{\lambda}\norm{a_2}_*+\norm{\varrho_{0}}_{\mathcal{C}^1(M)} \norm{a_2}_*+\lambda^2 \norm{\varrho_{0}}_{\mathcal{C}^1(M)}
   \norm{a_2}_*}\\ \nonumber
&\leq  C\para{\lambda \norm{\varrho_{0}}_{\mathcal{C}^1(M)}+\frac{1}{\lambda^2}}\norm{a_2}_*.
\end{align*}
Moreover, we have
\begin{align*}
\norm{\nabla v_\lambda (t,\cdot)}_{L^2(\M)}
&\leq C\eta \bigg\{\int_0^T\para{\lambda\norm{k_0(2\lambda t,\cdot)}_{L^2(\M)}+\norm{\p_t k_0(2\lambda t,\cdot)}_{L^2(\M)}} \, \dd t\cr
&\quad+\int_0^T\para{\lambda^2\norm{k_1(2\lambda t,\cdot)}_{L^2(\M)}+\lambda \norm{\p_t k_1(2\lambda t,\cdot)}_{L^2(\M)}} \, \dd t \cr
&\quad+ \int_0^T\para{\lambda^2\norm{k_2(2\lambda t,\cdot)}_{L^2(\M)}+\lambda \norm{\p_t k_2(2\lambda t,\cdot)}_{L^2(\M)}} \, \dd t  \bigg\} \cr
&\quad
\begin{aligned}
    +\frac{C}{\eta}\int_0^T\bigg(\frac{1}{\lambda}\norm{k_0(2\lambda t,\cdot)}_{L^2(\M)} &+ \norm{k_1(2\lambda t,\cdot)}_{L^2(\M)} \\
    &+\norm{k_2(2\lambda t,\cdot)}_{L^2(\M)}\bigg) \, \dd t.
\end{aligned}
\end{align*}
Hence, we obtain the following estimate
\begin{align}\label{6.25}
\|\nabla v_\lambda(t,\cdot)&\|_{L^2(\M)} \\ \nonumber
&\leq C\eta \bigg\{\int_\R\para{\norm{k_0(s,\cdot)}_{L^2(\M)}+\frac{1}{\lambda}\norm{\p_t k_0(s,\cdot)}_{L^2(\M)}} \, \dd t\cr
&\quad+\int_\R\para{\lambda\norm{k_1(s,\cdot)}_{L^2(\M)}+ \norm{\p_t k_1(s,\cdot)}_{L^2(\M)}} \, \dd t\cr
&\quad+ \int_\R\para{\lambda\norm{k_2(s,\cdot)}_{L^2(\M)}+\norm{\p_t k_2(s,\cdot)}_{L^2(\M)}} \, \dd t \bigg\} \cr
&\quad
\begin{aligned}
     +\frac{C}{\eta}\int_\R\bigg(\frac{1}{\lambda^2}\norm{k_0(s,\cdot)}_{L^2(\M)} &+ \frac{1}{\lambda}\norm{k_1(s,\cdot)}_{L^2(\M)} \\
     &+ \frac{1}{\lambda}\norm{k_2(s,\cdot)}_{L^2(\M)}\bigg) \,\dd t
\end{aligned}\cr
&\leq C\eta\para{1+\lambda^3\norm{\varrho_{0}}_{\mathcal{C}^1(M)}}\norm{a_2}_*+\frac{C}{\eta}
\para{\frac{1}{\lambda^2}+\lambda \norm{\varrho_{0}}_{\mathcal{C}^1(M)}}\norm{a_2}_*.
\end{align}
Now choosing $\eta=\lambda^{-1}$, we obtain
\begin{equation}\label{6.26}
\norm{\nabla v_\lambda(t,\cdot)}_{L^2(\M)} \leq C\para{\frac{1}{\lambda}+\lambda^2 \norm{\varrho_{0}}_{\mathcal{C}^1(M)}}\norm{a_2}_*.
\end{equation}
Finally, if we study the equation satisfied by $\p_tv$, we also find
\begin{equation}\label{6.27}
\norm{\p_t v(t,\cdot)}_{L^2(\M)}\leq C\para{1+\lambda^3  \norm{\varrho_{0}}_{\mathcal{C}^1(M)}} \norm{a_2}_*.
\end{equation}
This ends the proof of Lemma \ref{L6.2}.
\end{Demo}
\begin{Lemm}\label{L6.3}
There exists $C>0$ such that for any $a_1$, $a_2\in H^1(\R,H^2(\M))$ satisfying the transport equation \eqref{6.9} with \eqref{4.3}, the
following estimate holds true
\begin{multline}\label{6.28}
\lambda\abs{\int_0^T\!\!\!\int_\M \varrho(x)(a_1\overline{a}_2)(2\lambda t,x)\dv \dd t} \\
\leq C \bigg\{\norm{\varrho_0}_{\mathcal{C}^1(\M)} \para{\lambda^{-1}+\lambda^2 \norm{\varrho_0}_{\mathcal{C}^1(\M)}}
+\lambda\norm{\Lambda_{\g}-\Lambda_{\cg}}\bigg\}\norm{a_1}_*\norm{a_2}_*
\end{multline}
for any sufficiently large $\lambda$.
\end{Lemm}
\begin{Demo}{}
Following Lemma \ref{L6.2} let $u_2$ be a solution to the problem $(i\p_t+\Delta_{\cg})u=0$ of the form
$$
\overline{u}_2(t,x)=\frac{1}{\lambda}\overline{a}_2(2\lambda t,x)e^{-i\lambda(\psi_1-\lambda t)}+\overline{a}_3
(2\lambda t,x;\lambda)e^{-i\lambda(\psi_2-\lambda t)}+\overline{v}_{2,\lambda}(t,x)
$$
where $v_{2,\lambda}$ satisfies (\ref{6.15}) and $a_3$ satisfies (\ref{6.11}).
Thanks to Lemma \ref{L4.1} let $u_1$ be a solution of $(i\p_t+\Delta_\g)u=0$ of the form
$$
u_1(t,x)=a_1(2\lambda t,x)e^{i\lambda(\psi_1-\lambda t)}+v_{1,\lambda}(t,x).
$$
where $v_{1,\lambda}$ satisfies (\ref{4.7}). Then we have
\begin{align}\label{6.29}
\p_t\overline{u}_2(t,x)&= 2\p_t\overline{a}_2(2\lambda t,x)e^{-i\lambda(\psi_1-\lambda t)}+i\lambda
\overline{a}_2(2\lambda t,x)e^{-i\lambda(\psi_1-\lambda t)}\cr
&\quad +2\lambda\p_t\overline{a}_3(2\lambda t,x;\lambda)e^{-i\lambda(\psi_2-\lambda t)}+i\lambda^2\overline{a}_3
(2\lambda t,x,\lambda)e^{-i\lambda(\psi_2-\lambda t)}\cr
&\quad+\p_t\overline{v}_{2,\lambda}(t,x).
\end{align}
Let us compute the first term in the right hand side of (\ref{6.5}). We have
\begin{equation}\label{6.31}
\int_0^T\!\!\!\int_\M\varrho_1u_1\p_t\overline{u}_2\dv \, \dd t=
i\lambda\int_0^T\!\!\!\int_\M \varrho_1 (a_1\overline{a}_2)(2\lambda t,x)\dv \, \dd t+\mathcal{J}_1(\lambda) + \mathcal{J}_2(\lambda)
\end{equation}
with
\begin{align*}
\mathcal{J}_1(\lambda)&=2\int_0^T\!\!\!\int_\M \varrho_1\para{a_1\p_t\overline{a}_2}(2\lambda t,x)\dv \, \dd t \displaybreak[1] \\
&\quad+2\int_0^T\!\!\!\int_\M \varrho_1 v_{1,\lambda}\p_t\overline{a}_2(2\lambda t,x)e^{-i\lambda(\psi_1-\lambda t)}\dv \, \dd t \displaybreak[1] \\
&\quad+i\lambda\int_0^T\!\!\!\int_\M \varrho_1v_{1,\lambda}\overline{a}_2(2\lambda t,x)e^{-i\lambda(\psi_1-\lambda t)}\dv \, \dd t
\end{align*}
and with
\begin{align*}
\mathcal{J}_2(\lambda)&=
+2\lambda\int_0^T\!\!\!\int_\M \varrho_1\para{a_1\p_t\overline{a}_3}(2\lambda t,x)e^{i\lambda(\psi_1-\psi_2)}\dv \, \dd t \displaybreak[1]\\
&\quad+i\lambda^2\int_0^T\!\!\!\int_\M \varrho_1\para{a_1\overline{a}_3}(2\lambda t,x)e^{i\lambda(\psi_1-\psi_2)}\dv \, \dd t \displaybreak[1] \\
&\quad+\int_0^T\!\!\!\int_\M \varrho_1a_1(2\lambda t,x)\p_t\overline{v}_{2,\lambda}(t,x)\dv \, \dd t \displaybreak[1] \\
&\quad+2\lambda\int_0^T\!\!\!\int_\M \varrho_1v_{1,\lambda}\p_t\overline{a}_3(2\lambda t,x)e^{-i\lambda(\psi_2-\lambda t)}\dv \, \dd t \displaybreak[1]\\
&\quad+i\lambda^2\int_0^T\!\!\!\int_\M \varrho_1v_{1,\lambda}\overline{a}_3(2\lambda t,x)e^{-i\lambda(\psi_2-\lambda t)}\dv \, \dd t\\
&\quad+\int_0^T\!\!\!\int_\M \varrho_1v_{1,\lambda}\p_t\overline{v}_{2,\lambda}(t,x)\dv \,  \dd t.
\end{align*}
>From (\ref{6.11}), (\ref{6.15}) and (\ref{4.7}) we have the estimates
\begin{align}\label{6.32}
\abs{\mathcal{J}_1(\lambda)}&\leq C\norm{\varrho_0}_{\mathcal{C}^1(\M)}\lambda^{-1}\norm{a_2}_*\norm{a_1}_* \nonumber \\
\abs{\mathcal{J}_2(\lambda)}&\leq C\norm{\varrho_0}_{\mathcal{C}^1(\M)}\para{\lambda^{-1}
+\lambda^2\norm{\varrho_0}_{\mathcal{C}^1(\M)}}\norm{a_2}_*\norm{a_1}_*.
\end{align}
On the other hand, we have
\begin{align*}
\nabla_\g u_1&=(\nabla_\g a_1)(2\lambda t,x)e^{i\lambda(\psi_1-\lambda t)}+
i\lambda(\nabla_\g \psi_1)a_1(2\lambda t,x)e^{i\lambda(\psi_1-\lambda t)}+\nabla_\g v_{1,\lambda} \\
\nabla_\g \overline{u}_2&=\frac{1}{\lambda}(\nabla_\g\overline{a}_2)(2\lambda t,x)e^{-i\lambda(\psi_1-\lambda t)}-
i\overline{a}_2(2\lambda t,x)\nabla_\g\psi_1e^{-i\lambda(\psi_1-\lambda t)}\cr
&\quad-i\lambda\overline{a}_3(2\lambda t,x)
\nabla_{\g}\psi_2 e^{-i\lambda(\psi_2-\lambda t)}+\nabla_\g\overline{a}_3(2\lambda t,x)
e^{-i\lambda(\psi_2-\lambda t)}+\nabla_\g\overline{v}_{2,\lambda}
\end{align*}
and the second term in the right-hand side of (\ref{6.5}) becomes
\begin{multline}\label{6.35}
\int_0^T\!\!\!\int_\M \varrho_2(x)\seq{\nabla_\g u_1(t,x),\nabla_\g \overline{u}_2(t,x)}_\g \dv \, \dd t \\
=\lambda\int_0^T\!\!\!\int_\M \varrho_2(x)(a_1\overline{a}_2)(2\lambda t,x)\dv \, \dd t+\mathcal{J}_3(\lambda)+\mathcal{J}_4(\lambda)
\end{multline}
with
\begin{align*}
\mathcal{J}_3(\lambda) &= \frac{1}{\lambda}\int_0^T\!\!\!\int_\M\varrho_2(x)\seq{\nabla_\g a_1(2\lambda t,x),\nabla_\g\overline{a}_2(2\lambda t,x)}_\g\dv
\, \dd t \\
&\quad-i\int_0^T\!\!\!\int_\M\varrho_2(x)\overline{a}_2(2\lambda t,x)\seq{\nabla_\g a_1(2\lambda t,x),\nabla_\g\psi_1(x)}_\g\dv \,  \dd t \displaybreak[1] \\
&\quad+i\int_0^T\!\!\!\int_\M \varrho_2(x) a_1(2\lambda t,x)\seq{\nabla_\g\overline{a}_2(2\lambda t,x), \nabla_\g\psi_1(x)}_\g\dv \, \dd t\\
&\quad+\frac{1}{\lambda}\int_0^T\!\!\!\int_\M \varrho_2(x) e^{-i\lambda(\psi_1-\lambda t)}\seq{\nabla_\g\overline{a}_2(2\lambda t,x),
\nabla_\g v_{1,\lambda}(t,x)}_\g\dv \, \dd t\\
&\quad-i\int_0^T\!\!\!\int_\M \varrho_2(x) \overline{a}_2(2\lambda t,x)e^{-i\lambda(\psi_1-\lambda t)}\seq{\nabla_\g v_{1,\lambda}(t,x),
\nabla_\g\psi_1(x)}_\g\dv \, \dd t
\end{align*}
and with
\begin{align*}
\mathcal{J}_4(\lambda) &= -i\lambda\int_0^T\!\!\!\int_\M \varrho_2(x)\overline{a}_3(2\lambda t,x)e^{i\lambda(\psi_1-\psi_2)}\seq{\nabla_\g
a_1(2\lambda t,x),
\nabla_\g\psi_2(x)}_\g\dv \, \dd t \displaybreak[1] \\
&\quad+\int_0^T\!\!\!\int_\M\varrho_2(x) e^{i\lambda(\psi_1-\psi_2)}\seq{\nabla_\g a_1(2\lambda t,x),\nabla_\g\overline{a}_3(2\lambda t,x)}_\g\dv \,
\dd t
\displaybreak[1]\\
&\quad+\int_0^T\!\!\!\int_\M \varrho_2(x) e^{i\lambda(\psi_1-\lambda t)}\seq{\nabla_\g a_1(2\lambda t,x),\nabla_\g\overline{v}_{2,\lambda}(t,x)}_\g\dv
\, \dd t \displaybreak[1]\\
&\quad+\lambda^2\int_0^T\!\!\!\int_\M \varrho_2(x) (a_1\overline{a}_3)(2\lambda t,x)e^{i\lambda(\psi_1-\psi_2)}\seq{\nabla_\g\psi_1(x),
\nabla_\g\psi_2(x)}_\g\dv \, \dd t\\
&\quad+i\lambda\int_0^T\!\!\!\int_\M \varrho_2(x) a_1(2\lambda t,x)e^{i\lambda(\psi_1-\psi_2)}\seq{\nabla_\g\overline{a}_3(2\lambda t,x),
\nabla_\g\psi_1(x)}_\g\dv \, \dd t\\
&\quad+i\lambda\int_0^T\!\!\!\int_\M \varrho_2(x) a_1(2\lambda t,x)e^{i\lambda(\psi_1-\lambda t)}\seq{\nabla_\g\psi_1(x),\nabla_\g\overline{v}_2(t,x)}_\g
\dv \, \dd t\\
&\quad-i\int_0^T\!\!\!\int_\M \varrho_2 \overline{a}_2(2\lambda t,x))e^{-i\lambda(\psi_1-\lambda t)}\seq{\nabla_\g v_{1,\lambda}(t,x),
\nabla_\g\psi_1(x)}_\g\dv \, \dd t\\
&\quad-i\lambda\int_0^T\!\!\!\int_\M \varrho_2\overline{a}_3(2\lambda t,x)e^{-i\lambda(\psi_2-\lambda t)}\seq{\nabla_\g v_{1,\lambda}(t,x),
\nabla_g\psi_2}_\g\dv \, \dd t\\
&\quad+\int_0^T\!\!\!\int_\M \varrho_2 e^{-i\lambda(\psi_2-\lambda t)}\seq{\nabla_\g v_{1,\lambda},\nabla_\g \overline{a}_3(2\lambda t,x)}_\g\dv \,
\dd t\\
&\quad+\int_0^T\!\!\!\int_\M \varrho_2 \seq{\nabla_\g v_{1,\lambda}, \nabla_\g \overline{v}_{2,\lambda}}_\g\dv \, \dd t.
\end{align*}
>From (\ref{6.11}), (\ref{6.15}) and (\ref{4.7}), we have
\begin{align}\label{6.36}
\abs{\mathcal{J}_3(\lambda)}&\leq \norm{\varrho_0}_{\mathcal{C}^1(\M)}\lambda^{-1} \norm{a_2}_*\norm{a_1}_* \displaybreak[1] \\ \nonumber
\abs{\mathcal{J}_4(\lambda)}&\leq \norm{\varrho_0}_{\mathcal{C}^1(\M)}\para{\lambda^{-1}
+\lambda^2 \norm{\varrho_0}_{\mathcal{C}^1(\M)}}\norm{a_2}_*\norm{a_1}_*.
\end{align}
Taking into account (\ref{6.5}), (\ref{6.31}) and (\ref{6.35}), we deduce that
\begin{multline}\label{6.37}
\int_0^T\!\!\!\int_{\p\M}\para{\Lambda_{\g}-\Lambda_{\cg}}f_1\overline{f}_2\ds \,  \dd t=
\lambda\int_0^T\!\!\!\int_\M \varrho(x)(a_1\overline{a}_2)(2\lambda t,x)\dv \, \dd t \\ +
\mathcal{J}_1(\lambda)+\mathcal{J}_2(\lambda)+\mathcal{J}_3(\lambda)+\mathcal{J}_4(\lambda).
\end{multline}
In view of (\ref{6.32}) and (\ref{6.36}), we obtain
\begin{multline*}
\lambda\abs{\int_0^T\!\!\!\int_\M \varrho(x)(a_1\overline{a}_2)(2\lambda t,x)\dv \, \dd t} \\
\leq C\Big\{\norm{\varrho_0}_{\mathcal{C}^1(\M)}\para{\lambda^{-1}+\lambda^2\norm{\varrho_0}_{\mathcal{C}^1(\M)}}
+\lambda\norm{\Lambda_{\g}-\Lambda_{\cg}}\Big\}\norm{a_1}_{*}\norm{a_2}_{*}.
\end{multline*}
This completes the proof.
\end{Demo}
\subsection{Stability estimate for the geodesic ray transform}
\begin{Lemm}\label{L6.4}
There exists $C>0$ such that for any $b\in H^2(\p_+S\M_{1})$ the following estimate
\begin{multline}\label{6.39}
\abs{\int_{\p_+S\M_1}\I(\varrho)(y,\theta)b(y,\theta)\seq{\theta,\nu(y)}\dss (y,\theta)}\\
\leq C\Big(\big(\lambda^{-1}+\lambda^{2}\norm{\varrho_0}_{\mathcal{C}^1(\M)}\big)\norm{\varrho_0}_{\mathcal{C}^1(\M)}+\lambda
\norm{\Lambda_{\g}-\Lambda_{\cg}}\Big)\norm{b}_{H^2(\p_+S\M_{1})}
\end{multline}
holds for any $y\in\p\M_1$.
\end{Lemm}
\begin{Demo}{}
Following (\ref{4.21}), we take two solutions of the form
\begin{align*}
\widetilde{a}_1(t,r,\theta)&=\alpha^{-1/4}\phi(t-r)b(y,\theta), \\
\widetilde{a}_2(t,r,\theta)&=\alpha^{-1/4}\phi(t-r)\mu(y,\theta).
\end{align*}
Now we change variable in (\ref{6.28}), $x=\exp_{y}(r\theta)$, $r>0$ and
$\theta\in S_{y}\M_1$. Then
\begin{align*}
\int_0^T\!\!\int_\M &\varrho a_1(2\lambda t,x)a_2(2\lambda t,x)\dv\,\dd t\\
&=\int_0^T\!\!\int_{S_{y}\M_1}\!\int_0^{\tau_+(y,\theta)}\widetilde{\varrho}(r,\theta)\widetilde{a}_1(2\lambda t,r,\theta)
\widetilde{a}_2(2\lambda t,r,\theta)\alpha^{1/2} \,\dd r \, \dd \omega_y(\theta) \, \dd t \displaybreak[1] \\
&=\int_0^T\!\!\int_{S_{y}\M_1}\!\int_0^{\tau_+(y,\theta)}\widetilde{\varrho}(r,\theta)\phi^2(2\lambda t-r)b(y,\theta) \mu(y,\theta)
\, \dd r \, \dd \omega_y(\theta) \, \dd t \displaybreak[1] \\
&=\frac{1}{2\lambda}\int_0^{2\lambda T}\int_{S_{y}\M_1}\!\int_0^{\tau_+(y,\theta)}\widetilde{\varrho}(r,\theta)
\phi^2( t-r)b(y,\theta)\mu(y,\theta) \, \dd r \, \dd \omega_y(\theta) \, \dd t.
\end{align*}
We conclude that
\begin{multline}\label{6.41}
\abs{\int_0^T\!\!\int_{S_{y}\M_1}\!\int_0^{\tau_+(y,\theta)}\widetilde{\varrho}(r,\theta)\phi^2( t-r)b(y,\theta)\mu(y,\theta)
\, \dd r \,\dd \omega_y(\theta) \, \dd t} \\ \leq
C\Big(\big(\lambda^{-1}+\lambda^2\norm{\varrho_0}_{\mathcal{C}^1(\M)}\big)\norm{\varrho_0}_{\mathcal{C}^1(\M)} \\ +\lambda
\norm{\Lambda_{\g,\,q_1}-\Lambda_{\g,\,q_2}}\Big)\norm{b(y,\cdot)}_{H^2(S^+_y\M_{1})}
\end{multline}
where $S^+_y\M_{1}=\{\theta \in S_{y}\M: \langle \theta,\nu \rangle_{\g}<0\}$.
Given the support properties of the function $\phi$, the left-hand side of the inequality reads in fact
\begin{multline*}
\int_{-\infty}^{\infty}\int_{S_{y}\M_1}\!\int_0^{\tau_+(y,\theta)}\widetilde{\varrho}(r,\theta)\phi^2( t-r)b(y,\theta)\mu(y,\theta)
\, \dd r \, \dd \omega_y(\theta) \, \dd t \\ = \bigg(\int_{-\infty}^{\infty} \phi^2(t) \, \dd t\bigg) \times
\int_{S_{y}\M_1}\!\int_0^{\tau_+(y,\theta)}\widetilde{\varrho}(r,\theta)b(y,\theta)\mu(y,\theta)
\, \dd r \, \dd \omega_y(\theta).
\end{multline*}
Integrating with respect to $y\in \p \M_1$ in \eqref{6.41} we obtain
\begin{multline*}
\abs{\int_{\p_+S\M_1}\I(\varrho)(y,\theta)b(y,\theta)\seq{\theta,\nu(y)}\dss (y,\theta)}\cr
\leq \Big(\big(\lambda^{-1}+\lambda^{2}\norm{\varrho_0}_{\mathcal{C}^1(\M)}\big)\norm{\varrho_0}_{\mathcal{C}^1(\M)}+\lambda
\norm{\Lambda_{\g}-\Lambda_{\cg}}\Big)\norm{b}_{H^2(\p_+S\M_{1})}.
\end{multline*}
This completes the proof of the lemma.
\end{Demo}
\subsection{End of the proof of Theorem \ref{Th2}}
Let us now prove Theorem \ref{Th2}. We choose
$$
b(y,\theta)=\I\para{\I^*\I(q)}(y,\theta)
$$
and obtain using Lemma \ref{L6.4} and \eqref{3.4}
\begin{multline*}
\norm{\I^*\I(\varrho)}^2_{L^2(\M_1)} \leq C\Big(\big(\lambda^{-1}+\lambda^{2}\norm{\varrho_0}_{\mathcal{C}^1(\M)}\big)
\norm{\varrho_0}_{\mathcal{C}^1(\M)} \\ 
+\lambda\norm{\Lambda_{\g}-\Lambda_{\cg}}\Big)\norm{\I^*\I(\varrho)}_{H^2(\M_1)}.
\end{multline*}
By interpolation we have
\begin{align*}
\norm{\I^*\I(\varrho)}^2_{H^1(\M_1)}&\leq C \norm{\I^*\I(\varrho)}_{L^2(\M_1)}\norm{\I^*\I(\varrho)}_{H^2(\M_1)}\cr
&\begin{aligned}
    \leq C \Big(\big(\lambda^{-1}+\lambda^{2}\norm{\varrho_0}_{\mathcal{C}^1(\M)}\big)\norm{\varrho_0}_{\mathcal{C}^1(\M)}& \cr
    +\lambda \norm{\Lambda_{\g}-\Lambda_{\cg}}&\Big)^{1/2}\norm{\I^*\I(\varrho)}_{H^2(\M)}^{3/2}.
\end{aligned}
\end{align*}
We use \eqref{3.7} and \eqref{3.8} to deduce
\begin{align*}
\|\varrho\|^2_{L^2(\M)} &\lesssim \Big(\big(\lambda^{-1}+\lambda^{2}\norm{\varrho_0}_{\mathcal{C}^1(\M)}\big)\norm{\varrho_0}_{\mathcal{C}^1(\M)}
    +\lambda \norm{\Lambda_{\g}-\Lambda_{\cg}}\Big)^{\frac{1}{2}}\norm{\varrho}_{H^1(\M)}^{\frac{3}{2}} \\
&\lesssim \big(\lambda^{-1}+\lambda^{2}\norm{\varrho_0}_{\mathcal{C}^1(\M)}\big)^{\frac{1}{2}}\norm{\varrho_0}^2_{\mathcal{C}^1(\M)}
    +\lambda^{\frac{1}{2}}\norm{\varrho_{0}}_{\mathcal{C}^1(\M)}^{\frac{3}{2}} \norm{\Lambda_{\g}-\Lambda_{\cg}}^{\frac{1}{2}}.
\end{align*}
Minimizing $\lambda^{-1}+\lambda^{2}\norm{\varrho_0}_{\mathcal{C}^1(\M)}$ in $\lambda$,  we get
\begin{align*}
\norm{\varrho}^2_{L^2(\M)} &\lesssim \norm{\varrho_0}_{\mathcal{C}^1(\M)}^{13/6}
    +  \norm{\varrho_0}_{\mathcal{C}^1(\M)} \norm{\Lambda_{\g}-\Lambda_{\cg}}^{1/2} \\
&\lesssim \varepsilon^{1/12} \norm{\varrho_0}_{\mathcal{C}^1(\M)}^{25/12}
    +  \varepsilon \norm{\Lambda_{\g}-\Lambda_{\cg}}^{1/2}.
\end{align*}
Since
\begin{align*}
\norm{\varrho_0}_{\mathcal{C}^1(\M)}&\lesssim \norm{\varrho_0}_{H^{n/2+1+\epsilon}(\M)} \\  &\lesssim \norm{\varrho_0}_{L^2(\M)}^{24/25}
\norm{\varrho_0}_{H^s(\M)}^{1/25} \lesssim \norm{\varrho_0}_{L^2(\M)}^{24/25}
\end{align*}
we therefore obtain
$$
\norm{\varrho}^2_{L^2(\M)}\lesssim \varepsilon^{1/12} \norm{\varrho_0}_{L^2(\M)}^{2}
    +   \norm{\Lambda_{\g}-\Lambda_{\cg}}^{1/2}.
$$
Taking $\varepsilon>0$ small enough we conclude and obtain (\ref{1.14}).

\subsection*{Acknowledgements} Part of this work was done while David DSF was visiting the Facult\'e des Sciences de Bizerte.
He gratefully acknowledges  the hospitality of the Facult\'e des Sciences.

\end{document}